\documentclass{amsart}
\usepackage{amsmath}
\usepackage{amssymb}

\usepackage{mathtools}
\usepackage{color,titletoc}
\usepackage[shortlabels]{enumitem}
\usepackage{mathrsfs,accents}

\usepackage{natbib}
\usepackage[colorlinks=true,breaklinks=true,bookmarks=true,urlcolor=blue,
citecolor=blue,linkcolor=blue,bookmarksopen=false,draft=false]{hyperref}
\usepackage[justification=centering]{caption}
\usepackage{mathrsfs}
\usepackage[normalem]{ulem}

\usepackage{worldflags}
\usepackage{amsfonts}
\usepackage{mathtools}
\usepackage{xcolor}
\usepackage{algorithm}
\usepackage{algorithmic}
\usepackage{latexsym}
\usepackage{booktabs}
\usepackage{textcomp}
\usepackage{subfigure}
\usepackage{enumitem}
\usepackage{amsmath}
\usepackage{longtable}

\bibpunct[, ]{[}{]}{,}{n}{}{,}%

\newcommand{\re}{\mathbb{R}}

\newcommand{\CRed}{\textcolor{black}}

\newcommand{\N}{\mathbb{N}}
\newcommand{\Z}{\mathbb{Z}}

\def\af{\alpha}

\newcommand{\st}{\mathrm{s.t.}}
\newcommand{\reff}[1]{\eqref{#1}}
\newcommand{\lmd}{\lambda}

\newcommand{\dt}{\delta}

\newcommand{\eps}{\epsilon}

\newcommand{\mt}[1]{\mathtt{#1}}
\newcommand{\QM}[1]{\operatorname{QM}[#1]}
\newcommand{\ideal}[1]{\operatorname{Ideal}[#1]}
\newcommand{\qmod}[1]{\operatorname{QM}[#1]}
\renewcommand{\top}{T}

\newcommand{\nn}{\nonumber}

\newcommand{\sig}{\sigma}

\newcommand{\supp}[1]{\mathrm{supp}(#1)}

\DeclareMathOperator{\relint}{relint}
\DeclareMathOperator{\rank}{rank}

\def\ben{\begin{enumerate} }%
	\def\een{\end{enumerate}}

\def\bem{\begin{pmatrix}}
	\def\eem{\end{pmatrix}}

\def\beq{\begin{equation}}
	\def\eeq{\end{equation}}

\renewcommand{\subset}{\subseteq}

\renewcommand{\emptyset}{\varnothing}

\newcommand{\T}{T}

\def\al{\alpha}
\def\be{\beta}

\def\epsilon{\varepsilon}
\def\eps{\epsilon}

\def\bex{\begin{example}}
	\def\eex{\end{example}}

\def\R{ {\mathbb{R}} }
\def\C{ {\mathbb{C}} }

\def\N{ {\mathbb{N}} }

\def\la{\lambda}

\newcommand{\x}{{\tt x}}
\newcommand{\y}{{\tt y}}
\newcommand\cx{[\x]}
\newcommand\cy{[\y]}
\newcommand\cxy{[\x,\y]}
\newcommand\rx{\R\cx}
\newcommand\ry{\R\cy}
\newcommand\rxy{\R\cxy}

\newcommand{\dd}{{\mathrm d}}

\newcommand\sosx{\Sigma^2[\x]}

\newcommand\sosy{\Sigma^2[\y]}
\newcommand\sosxy{\Sigma^2[\x,\y]}

\def\beq{\begin{equation} }%
\def\eeq{\end{equation}}

\def\ben{\begin{enumerate} }%
	\def\een{\end{enumerate}}

\def\bem{\begin{pmatrix}}
	\def\eem{\end{pmatrix}}

\def\beq{\begin{equation}}
	\def\eeq{\end{equation}}

\renewcommand{\subset}{\subseteq}

\renewcommand{\emptyset}{\varnothing}

\def\bex{\begin{example}}
	\def\eex{\end{example}}

\numberwithin{equation}{section}

\newtheorem{thm}{Theorem}[section]
\newtheorem{theorem}[thm]{Theorem}
\newtheorem{corollary}[thm]{Corollary}

\newtheorem{prop}[thm]{Proposition}
\theoremstyle{definition}

\newtheorem{assum}[thm]{Assumption}
\newtheorem{remark}[thm]{Remark}

\newtheorem{example}[thm]{Example}

\numberwithin{equation}{section}

\begin{document}

\title[Positivstellens\"atze and Moment problems with Quantifiers]
{Positivstellens\"atze and Moment problems\\ with Universal Quantifiers}

\author{Xiaomeng Hu}
\address{Xiaomeng Hu,
Department of Mathematics,
University of California San Diego,
9500 Gilman Drive, La Jolla, CA, USA, 92093}
\email{x8hu@ucsd.edu}

\author{Igor Klep}
\address{Igor Klep, Faculty of Mathematics and Physics,
University of Ljubljana \& Faculty of Mathematics, Science
and Information Technology, University of Primorska, Koper \& Institute of Mathematics,
    Physics and Mechanics, Ljubljana, Slovenia}
\email{igor.klep@fmf.uni-lj.si}

\author{Jiawang Nie}
\address{Jiawang Nie, Department of Mathematics,
University of California San Diego,
9500 Gilman Drive, La Jolla, CA, USA, 92093}
\email{njw@math.ucsd.edu}

\date{\today}


\begin{abstract}
This paper studies Positivstellens\"atze and moment problems for sets $K$
that are given by universal quantifiers.
Let $Q\subset \re^m$ be a closed set and
let $g = (g_1,\ldots,g_s)$ be a tuple of polynomials in
two vector variables {$\x=(\x_1,\ldots,\x_n)$ and $\y = (\y_1,\ldots,\y_m)$.}
Then $K$ is described as the set of all points $x\in \re^n$ such that
each $g_j(x, y) \ge 0$ for all $y \in Q$.
Fix a {finite nonnegative Borel} measure $\nu$ on $\re^m$ with $\supp{\nu} = Q$,
and assume it satisfies the {multivariate} Carleman condition.

The first main result of the paper is a Positivstellensatz with universal quantifiers:
if a polynomial $f(\x)$
is positive on $K$, then $f(\x)$ belongs to the quadratic module
 $\QM{g,\nu}$ associated to $(g,\nu)$,
under the archimedeanness assumption on $\QM{g,\nu}$.
Here, $\QM{g,\nu}$ denotes the quadratic module of polynomials
in $\x$
that can be represented as
\[
\tau_0(\x) +  \int \tau_1(\x,y)g_1(\x, y)\, \dd \nu(y)   + \cdots  +
\int \tau_s(\x,y) g_s(\x, y)\, \dd \nu(y),
\]
where each $\tau_j$ is a sum of squares polynomial.

Second, necessary and sufficient conditions
for a full (or truncated) multisequence to admit a representing measure
supported in $K$ are given.
In particular, the classical flat extension
theorem of Curto and Fialkow is generalized to
truncated moment problems on such a set $K$.

Third, we present applications of the above Positivstellensatz and moment problems
in semi-infinite optimization, whose feasible sets are given by infinitely many
constraints with universal quantifiers. This results in a new hierarchy of
Moment-SOS relaxations. Its convergence is shown under some usual assumptions.
The quantifier set $Q$ is allowed to be non-semialgebraic, which makes it possible
to solve some optimization problems with non-semialgebraic constraints.  
\end{abstract}

\maketitle

\section{Introduction}
\label{sec:Intro}

Positivstellens\"atze and moment problems are
pillars of real algebraic geometry \cite{BCR98,Lau09,Scheid09}
and are of broad interest
in computational and applied  mathematics.
This paper concerns these two topics
when the constraining sets are
given by universal quantifiers.
Let $\x=(\x_1,\ldots,\x_n)$ and $\y=(\y_1,\ldots,\y_m)$ be tuples of
variables.
We are interested in {subsets} $K$ {of} $\re^n$
that are given by inequalities in $x$,
with $y$ as a universal quantifier {(see \cite{lasserre2015quantifier}).}
Let $Q \subseteq \re^m$ be a given closed set.
For a  tuple $g = (g_1,\ldots,g_s)$ of polynomials in $\rxy$,
consider the following set given by the universal quantifier $y$:\looseness=-1
\beq  \label{eq:K}
K=\{x\in\R^n :  g_1(x,y)\geq0,\ldots,g_s(x,y)\geq 0
\;\; \forall  y \in Q \}.
\eeq
When there is no universal quantifier $y$,
the set $K$ is a classical basic closed semialgebraic set.
By Tarski's transfer principle \cite{BCR98},
if the quantifier set {$Q$ is semialgebraic,
then $K$ is semialgebraic. A quantifier-free description for $K$
can be obtained by applying symbolic computations like cylindrical algebraic decompositions
(see \cite{BPR}). However, computing a quantifier-free description
is typically computationally expensive.
In this paper, the set $Q$ is allowed to be non-semialgebraic,
so $K$ may also be non-semialgebraic.
For instance, when $Q=\Z^m$ ($\Z$ denotes the set of integers),
the set $K$ is defined by countably many constraints.
\begin{itemize}

\item For $Q=\Z^1$ and $g(\x,\y) = (\x_1 +\y)^2 + \y^2 - \x_2^2 \ge 0$,
the set $K$ is given by
\[
x_1^2 - x_2^2 \ge 0, \quad 1 \ge \frac{x_2^2}{k^2} - \frac{(x_1+k)^2}{k^2}
\quad \text{for} \quad   k = \pm 1, \pm 2, \ldots.
\]

\item For $Q=\Z^1$ and $g(\x,\y) = \x_2 - 2 \y \x_1  + \y^2 \ge 0$, the set $K$ is
a convex polygon with infinitely many sides.

\item For $Q=\Z_+$ (the set of positive integers) and
$g(\x,\y) = 4\y^4 -1 - \y^2(2\y^2-1) \x_1^2 - \y^2(2\y^2+1) \x_2^2   \ge 0$,
the set $K$ is the intersection of infinitely many ellipses:
\[
\frac{x_1^2}{2+k^{-2}} + \frac{x_2^2}{2-k^{-2} }  \le 1 \quad
\text{for} \quad k = 1, 2, \ldots.
\]

\end{itemize}
   }
\looseness=-1

Positivstellens\"atze concern representations of polynomials
that are positive (or nonnegative) on a set $K$.
Equivalently, for a given polynomial $f \in \rx$, what is a
test or certificate for $f \ge 0$ on $K$?
When does such a certificate hold necessarily?
When $K$
has no universal quantifier $y$
(i.e., the polynomials $g_i$ in \eqref{eq:K} do not depend on $y$),
Positivstellens\"atze have been {extensively studied,}
see, e.g., the surveys and books
\cite{BCR98,MomentSOShierarchy,LasBk15,Lau09,MPO23,Scheid09}
or the following small sample of recent papers
\cite{CKS09, EP20, Fri21, GKKS15, LPR20, MNR23, PV99, Rie16, SS+, Scw03}
and the references therein.
For instance, consider the Putinar certificate\looseness=-1
\beq \label{rep:Put}
f \, = \, \sig_0 + \sig_1 g_1 + \cdots + \sig_s g_s,
\eeq
where all $\sig_i$ are sum-of-squares (SOS) polynomials in $\rx$.
Clearly, if $f$ has a representation of the form \reff{rep:Put},
then $f \ge 0$ on the set $K$. When the quadratic module of $g$
is archimedean, if $f>0$ on $K$,
then by Putinar's Positivstellensatz \cite{Putinar}
a representation of the form \reff{rep:Put} must hold.
A  representation more general than \reff{rep:Put}
is given by the Schm\"{u}dgen Positivstellensatz \cite{Schmu91},
which uses the preordering of $g$.
All these classical results assume that
$K$ is a basic closed semialgebraic set.
However, when $K$ depends on quantifiers as in \reff{eq:K},
there is little work on Positivstellens\"atze.
This is remedied in the present paper.

Closely related to Positivstellens\"atze are moment problems.
Let $\N$ denote the set of nonnegative integers
and $\N^n$ denote the set of nonnegative integer vectors of length $n$.
For a given multisequence $z = (z_\af)_{\af \in \N^n}$,
i.e., $z$ is a vector whose entries are labelled by
nonnegative integer vectors in $\N^n$,
the moment problem concerns the existence of
a {nonnegative} Borel measure\footnote{All our measures will be assumed finite.} $\mu$ on $\re^n$ such that
\beq\label{eq:MP}
z_\af \, = \, \int x^\af\ \dd\mu(x)  \quad  \forall \, \af \in \N^n.
\eeq
In the above, $x^\af \coloneqq x_1^{\af_1} \cdots x_n^{\af_n}$
for the multiindex
$\af = (\af_1, \ldots, \af_n)$.
The sequence $z$ is said to be a {\it moment sequence}
if such a Borel measure $\mu$ exists, and in this case
$\mu$ is called a representing measure for $z$.
We refer the reader to the surveys
\cite{Ber,Fia16}, books \cite{Akh,Scm},
or  papers \cite{BS16, BL20, CMN11, CGIK23, IK17, IKKM22, IKKM23, KW13, Net08, PS01}
and the references therein for more details about moment problems.

In many applications, the support of the measure $\mu$
is often required to be {contained} in a set $K$,
i.e., $\supp{\mu} \subseteq K$.
{Then $z$ is called a {\it $K$-moment sequence}
and $\mu$ is called a {\it $K$-representing measure} for $z$,}
if \eqref{eq:MP} holds for a Borel measure $\mu$ with $\supp{\mu} \subseteq K$.
When $K$ is described without quantifiers,
this is the classical $K$-moment problem
(see \cite{Fia16,Scm}).
However, there is little work on the $K$-moment problem
when $K$ depends on the quantifier $y$.
This is the second main topic of the present paper.\looseness=-1

Positivstellens\"{a}tze and moment problems with universal quantifiers
are useful for solving  semi-infinite optimization problems.
A typical problem of semi-infinite optimization is\looseness=-1
\beq \label{intro:SIP}
\left\{
\begin{array}{lll}
	\min\limits_{x\in X}& f(x)\\
	\st    &  g(x,y)\geq 0\quad\forall y\in Q.  \\
\end{array}
\right.
\eeq
Here, the constraining function $g$ depends on both $x$ and $y$
as is the case for the $g_j$ in \reff{eq:K},
and $X\subseteq\mathbb{R}^n$
is another given constraining set for $x$
that does not depend on the quantifier $y$.
The quantifier set $Q$ in \reff{intro:SIP}
need not be a basic closed semialgebraic set.
Solving this kind of semi-infinite optimization problem
is typically a highly challenging task.
However, Positivstellens\"{a}tze and moment problems with universal quantifiers
are powerful mathematical tools for solving them.
This is the third main topic of our paper.

\subsection*{Contributions}
The new contribution of this paper
is to solve the three above mentioned major {problems}.

Our first contribution is a Positivstellensatz
for sets $K$
defined with universal quantifiers as in \eqref{eq:K}.
{If a polynomial $f \in \re[\x]$ has the representation }
\beq\label{eq:rhsPsatz}
f(\x) = \sig(\x) +
\sum_{j=1}^s\int \tau_j(\x,y)g_j(\x,y)\, \dd \nu(y),
\eeq
where $\sig$ is an SOS polynomial in $\x$, $\tau_1, \ldots, \tau_s$
are SOS polynomials in $(\x,\y)$ and $\nu$ is a Borel measure on $\re^m$
such that $\supp{\nu} \subseteq Q$,
   {then we clearly have $f \ge 0$ on $K$.}
The Positivstellensatz {ensures}
that the reverse implication is \CRed{essentially} true.
The set of all polynomials in $\R\cx$ that can be written
as in \reff{eq:rhsPsatz}
is denoted by $\QM {g,\nu}$.
It is called the \textit{quadratic module generated by $g$ and $\nu$}.
Assume $\nu$ is a  Borel measure on $\re^m$ such that
$\supp{\nu} = Q$ and $\nu$
satisfies the Carleman condition %
\beq \label{fix:suppnu=Qintro}
\sum_{d=0}^\infty  \Big(\int y_j^{2d}\, \dd\nu(y) \Big)^{-\frac1{2d}}
= \infty \quad \text{for} \quad j=1, \ldots, m.
\eeq
We show in
Theorem \ref{thm:univPsatz}
that if $f >0$ on $K$
and $\QM {g,\nu}$ is archimedean
(i.e., $N-\x_1^2-\cdots -\x_n^2 \in \QM {g,\nu}$ for some positive integer $N$),
then $f \in \QM {g,\nu}$ must hold.
This is a generalization of Putinar's Positivstellensatz
to sets given by universal quantifiers.
Since the
truncations of the
quadratic module $\QM{g,\nu}$
for given degrees
can be represented by semidefinite programs (SDPs),
Theorem \ref{thm:univPsatz} gives rise to a Moment-SOS hierarchy of SDP relaxations
   to optimize a polynomial  over $K$.

Our second contribution is about {$K$-moment problems for sets $K$ defined
   by universal quantifiers {as in \eqref{eq:K}}.}
A key tool for studying moment problems is the {Riesz functional}.
A multisequence $z = (z_\al)_{\al \in \N^n}$
gives rise to the linear functional:
\beq   \nn
\mathscr{R}_z :  \rx\to\R, \quad \x^\al \mapsto z_\al .
\eeq
This is equivalent to
$
\mathscr{R}_z (f) = \sum_\af f_\af z_\af
$
for the polynomial $f = \sum_\af f_\af \x^\af$.
The functional $\mathscr{R}_z$ is called the \textit{Riesz functional} of $z$.
If $\mu$ is a representing measure for $z$, then
\[
\mathscr{R}_z(f) = \int f(x)\, \dd \mu(x)
\quad  \text{for all} \quad  f \in \rx.
\]
If in addition, $\supp{\mu} \subseteq K$, then
\beq\label{eq:Kpositive}
\mathscr{R}_z(f) \, \ge \, 0
\quad  \text{for all} \quad f \in \rx: \, f|_K \ge 0.
\eeq
We say the Riesz functional $\mathscr{R}_z$ is \textit{$K$-positive} if \reff{eq:Kpositive} holds.
{The Riesz functional $\mathscr{R}_z$
is simply said to be positive if it is $\re^m$-positive.}
The $K$-positivity is necessary
for $z$ to have a $K$-representing measure.
For a closed set $K$, being $K$-positive is also sufficient.
This is a classical result of M.~Riesz ($n=1$)
and Haviland ($n > 1$); see the works \cite{Akh,Ber,Fia16,Haviland,Riesz,Scm} for details.
When the quadratic module $\QM{g,\nu}$ is archimedean,
we show in Section \ref{sc:Mom} that $z$ is a $K$-moment sequence
if and only if $\mathscr{R}_z(f) \ge 0$ for all $f \in \QM{g,\nu}$.
Moreover, we also give concrete conditions for
$\mathscr{R}_z \ge 0$ on $\QM{g,\nu}$
in terms of moment and localizing matrices (see Theorem \ref{thm:shift}).

Our third contribution is on semi-infinite optimization.	
Suppose the constraining set\looseness=-1
\[
X \, = \, \{ x \in \re^n:\, c_{eq}(x) = 0,\, c_{in}(x) \ge 0\},
\]
for two tuples $c_{eq}, c_{in}$ of polynomials in $\x$.
{We refer to \reff{ideal:h} for the definition of the ideal $\ideal{c_{eq}}$ generated by $c_{eq}$
and refer to \reff{ideal:h} for the quadratic module $\qmod{c_{in}}$
generated by $c_{in}$.}
When $\qmod{g,\nu} + \ideal{c_{eq}} + \qmod{c_{in}}$ is archimedean,
we show in Section \ref{sc:SIP} that the semi-infinite optimization problem~\reff{intro:SIP}
is equivalent to
\begin{equation} 
\left\{ \begin{array}{rl}
	\min\limits_{{ z \in \re^{ \N^n } }  } & \mathscr{R}_z(f) \\
	\st  & \mathscr{R}_z \geq 0 \quad \mbox{on} \quad
	\qmod{g,\nu} + \ideal{c_{eq}} + \qmod{c_{in}},  \\
	& \mathscr{R}_z(1) = 1.
\end{array} \right.
\end{equation}
When the ideals and quadratic modules are
truncated by degrees, the above produces a hierarchy of
Moment-SOS type semidefinite programming relaxations.
We prove the convergence property for this hierarchy in Theorem \ref{thm:lass}.
Finally, we also discuss how to estimate moments of
the measure $\nu$ by sampling when the moments
are not known explicitly.
{We remark that the quantifier set $Q$ is allowed to be non-semialgebraic.
So this makes it possible to solve some semi-infinite optimization problems
with non-semialgebraic constraints.}

The paper is organized as follows.
Notation is fixed and some background on polynomial optimization and
moment problems is given in Section~\ref{sc:Pre}.
Positivstellens\"{a}tze, moment problems and semi-infinite optimization
for sets given by universal quantifiers
are respectively presented in Section~\ref{sc:Pos},
Section~\ref{sc:Mom}, and Section~\ref{sc:SIP}.
Some computational experiments are presented
in Section~\ref{sc:Num}.
Finally, in Section~\ref{sec:con}, we present our conclusions and engage
in a detailed discussion of our findings.\looseness=-1

\section{Preliminaries}
\label{sc:Pre}

\subsection{Notation}

The symbol $\rx=\mathbb{R}[\x_1,\ldots,\x_n]$ denotes the ring of  polynomials in
{$\x=(\x_1,\ldots,\x_n)$} with real coefficients.
The symbol $\re_{+}$ stands for the set of nonnegative real numbers.
For a symmetric matrix $W$,
$W\succeq 0$ means that $W$ is positive semidefinite. For a vector $u$,
$\|u\|$  denotes its standard Euclidean norm.
The notation $I_n$ denotes the $n\times n$ identity matrix.
The superscript $^T$ denotes the transpose of a matrix or vector.
The {symbol} $e$ denotes the vector of all ones, i.e.,
$e =(1, \ldots, 1)$.
We use $\otimes$ to denote the classical Kronecker product.

For {$\x=(\x_1,\ldots,\x_n)$} and
$\alpha \, \coloneqq \, (\alpha_1,\ldots,\alpha_n)\in\mathbb{N}^n$,
the notation {$\x^\alpha \coloneqq  \x^{\alpha_1}_1\cdots \x^{\alpha_n}_n$}
stands for the monomial of $\x$ with power $\alpha$.
We denote the power set
\[
\mathbb{N}^n_d =  \{\alpha\in\mathbb{N}^n \,\vert \,\alpha_1+\cdots+\alpha_n\leq d\} .
\]
Denote by $\mathbb{R}^{\mathbb{N}^n_{d}}$ the space of real vectors
that are labeled by $\alpha \in\mathbb{N}^n_{d}$.
For a positive integer {$d$}, the vector of all monomials in $\x$ of degrees at most $d$,
ordered with respect to the graded lexicographic ordering, is denoted as
\[
[\x]_{d} \, \coloneqq \,
\begin{pmatrix} 1 \ & \x_{1} \ & \cdots \  & \x_{n} \ & \x_{1}^{2} \ &
	\x_{1} \x_{2}  \   & \cdots\  & \x_{n}^{d}
    \end{pmatrix}^T.
\]

A polynomial $\sigma \in \mathbb{R}[\x]$ is said to be
a sum of squares (SOS) polynomial if
$\sigma=\sigma_{1}^{2}+\cdots+\sigma_{k}^{2}$
for some $\sigma_{1}, \ldots, \sigma_{k} \in \mathbb{R}[\x]$
   {and $k \in \N\setminus\{0\}$}.
The symbol $\Sigma^2[\x]$
denotes the cone of SOS polynomials in $\x$.
An interesting fact is that SOS polynomials
can be represented through semidefinite programming \cite{LasBk15,MPO23}.
Clearly, each SOS polynomial is nonnegative,
while not every nonnegative polynomial is SOS.
The approximation performance of SOS polynomials
is given in \cite{NieSOSbd}.
Moreover, SOS polynomials are also
very useful in tensor optimization
\cite{njwSTNN17,NieZhang18}.

For two sets $S, T \subseteq \mathbb{R}[\x]$,
their product and addition are defined as
\[
S \cdot T = \{ pq: \, p \in S, q \in T \}, \quad
S + T = \{ p+q: \, p \in S, q \in T \} .
\]
In particular, if $S = \{ p \}$ is a singleton, then we also use
\[
p \cdot T = \{ pq: \, q \in T \}, \quad
p + T = \{ p+q: \,  q \in T \} .
\]
A polynomial tuple $h=\left(h_{1}, \ldots, h_{m}\right)$ in $\mathbb{R}[\x]$
generates the ideal
\beq \label{ideal:h}
\ideal{h} \, \coloneqq \, h_{1} \cdot \mathbb{R}[\x]+\cdots+h_{m} \cdot \mathbb{R}[\x],
\eeq
which is the smallest ideal containing all $h_{i}$.
{For $k\in\mathbb{N}$ and $k\geq \deg(h)\coloneqq \max\{\deg(h_1),\ldots,\deg(h_m)\}$},
the $k$th truncation of $\ideal{h}$ is
\[
\ideal{h}_{k}  \coloneqq  h_{1} \cdot \mathbb{R}[x]_{k-\operatorname{deg}
	\left(h_{1}\right)}+\cdots+h_{m} \cdot
\mathbb{R}[x]_{k-\operatorname{deg}\left(h_{m}\right)}.
\]

A tuple $q=\left(q_{1}, \ldots, q_{t}\right)$
of polynomials in $\rx$ gives rise to the quadratic module
(let $q_0 \coloneqq 1$)
\beq \label{qmod:q}
\QM q \, \coloneqq \,
\Bigg \{ \sum_{i=0}^t \sigma_{i} q_{i} \,\Big \vert\,
\sigma_{i} \in \Sigma^2[\x]  \Bigg \}.
\eeq
For $k \in \mathbb{N}$ with $2 k \geq \operatorname{deg}(q)$,
the $k$th truncation of $\QM{q}$ is
\[
\QM{q}_{2k} \, \coloneqq \,
\Bigg \{\sum_{i=0}^{t} \sigma_{i} q_{i}\,\Big \vert\,
\sigma_{i} \in \Sigma^2[\x], \operatorname{deg}\left(\sigma_{i} q_{i}\right) \leq 2k
\Bigg \}.
\]
The quadratic module $\QM{q}$ is said to be {\it archimedean}
if there exists {an integer} $N>0$ such that
\[
N-\x_1^2 - \cdots - \x_n^2 \, \in \, \QM{q}.
\]
Quadratic modules are basic concepts in
polynomial optimization and moment problems.
We refer to
\cite{MomentSOShierarchy,LasBk15,Lau09,MPO23,Scheid09}
for recent work in this area.

Let $Q \subset \re^m$ be a closed set and
$\nu$ be a {nonnegative} Borel measure on $\re^m$ such that
$\supp{\nu} = Q$.  We let $L^2(\R^m,\nu)$
denote the Hilbert space of all $L^2$-integrable functions $\phi$ on $Q$, i.e.,
$\int \phi(y)^2\, \dd\nu(y) < \infty$.
The inner product on $L^2(\R^m,\nu)$ is given by
\[
\langle \phi, \psi \rangle_{L^2}  =
\int  \phi(y) \psi(y)\, \dd\nu(y), \quad
\phi,  \psi \in L^2(\R^m,\nu) .
\]
A linear functional $\ell$ on  $\re[\y]$, with $\y =(\y_1, \ldots, \y_m)$,
is said to satisfy the {\it multivariate Carleman condition} if
\beq \label{Carleman:fun}
	\sum_{d=0}^\infty  \Big(  \ell (\y_j^{2d} ) \Big)^{-\frac1{2d}}
	= \infty \quad \text{for} \quad j=1, \ldots, m.
\eeq

\section{Positivstellens\"atze with universal quantifiers}
\label{sc:Pos}

This section proves a Positivstellensatz
for polynomials $f$ positive on
a set $K$
given by a universal quantifier %
as in \reff{eq:K}.
Let $Q\subseteq\R^m$ be a given closed set.
We fix a {nonnegative} Borel measure $\nu$ on $\re^m$ satisfying the following assumption.

\begin{assum}\label{ass:measure}
The \CRed{nonnegative} Borel measure $\nu$ has the support $\supp{\nu} = Q$
and it satisfies the {\it multivariate Carleman condition}
\beq \label{fix:suppnu=Q}
	\sum_{d=0}^\infty  \Big(\int y_j^{2d}\, \dd\nu(y) \Big)^{-\frac1{2d}}
	= \infty \quad \text{for} \quad j=1, \ldots, m.
	\eeq
\end{assum}

A measure $\nu$ satisfying \reff{fix:suppnu=Q}
is known to be {\it determinate} (i.e., it is uniquely determined
by its moments $\int y^\alpha\, \dd\nu(y)$)
   {by Nussbaum's theorem {\cite{Nuss}},}
   and it is {\it strictly determinate} (i.e., $\ry$ is dense in $L^2(\R^m,\nu)$).
   See, e.g., \cite[Section 14.4]{Scm} for details and proofs.
It is interesting to remark that
the Carleman condition \reff{fix:suppnu=Q}
is {automatically satisfied} if $Q = \supp{\nu}$ is bounded.\looseness=-1

\subsection{Density of SOS polynomials }

In this subsection we prove the following strengthening of
the above-mentioned Nussbaum theorem:

\begin{prop}\label{prop:+approx}
Let $\nu$ be a {nonnegative Borel} measure satisfying Assumption~\ref{ass:measure}.
Then SOS polynomials are dense in the cone of
nonnegative functions in $L^2(\R^m,\nu)$.
\end{prop}
	
\proof{Proof.}
Suppose that the conclusion is not true.
Then there exists a nonnegative function $\phi\in L^2(\R^m,\nu)$
that is not in the $L^2$-closure of the convex cone $\sosy$.
By the Hahn-Banach separation theorem
(see \cite[Theorem III.3.4]{Bar02}),
there is a continuous linear functional $\ell:L^2(\R^m,\nu)\to\R$
satisfying
\beq  \label{eq:lpos}
\ell\big(\sosy\big)\subseteq\R_{+},\quad
\ell(\phi)<0.
\eeq
By adding a small multiple of the linear functional
$f\mapsto \int f\, {\rm d}\nu$ to $\ell$,
we can without loss of generality assume
there exists $\eps >0$ such that
\[
\ell(\sigma)\geq\eps>0 \quad \text{ for all } \sigma\in\sosy
\text{ with } \|\sigma\|_{L^2}=1.
\]
The Riesz representation theorem implies there is $h\in L^2(\R^m,\nu)$ such that
\[
\ell(f)  \,=\, \langle f,h\rangle_{L^2}  \,=\, \int f h \, {\rm d}\nu
\]
for all $f \in L^2(\R^m,\nu)$.
Since Assumption \ref{ass:measure} holds,
 $\ry$ is dense in $L^2(\R^m,\nu)$
(see, e.g., \cite[Theorem 14.2, Section 14.4]{Scm}).
Hence there is a sequence of polynomials $\{ p_n \}_{n=1}^\infty \subseteq \ry$
that converges to $h$ in the $L^2$-norm.
Applying the Cauchy-Schwartz inequality yields
\[
\Big| \langle f,h\rangle_{L^2} -
\langle f,p_n\rangle_{L^2} \Big| =
\Big| \langle f,h-p_n\rangle_{L^2} \Big|
\leq
\| f\|_{L^2}\, \|h-p_n\|_{L^2}.
\]
Hence, for $n$ large enough, the continuous linear functional
\beq\label{eq:ln}
\ell_n:f\mapsto\langle f,p_n\rangle_{L^2}
\eeq
also satisfies \eqref{eq:lpos}, i.e., $\ell_n$
is nonnegative on $\sosy$ while it is negative at $\phi$.

We now adapt the argument in \cite[Theorem 14.25]{Scm} to show
that $p_n\geq0$ on $\supp\nu$. The restriction $\ell_n:\ry\to\R$
is a positive linear functional and satisfies the multivariate
Carleman condition (see \cite[Corollary 14.22]{Scm}).
Hence, by Nussbaum’s theorem, $\ell_n$ is of the form\looseness=-1
\[
\ell_n(f) = \int f\,\dd\tau, \quad \forall f\in\ry
\]
for some  {nonnegative Borel} measure $\tau$ on $\R^m$. Set
\[
M_+ \coloneqq \{y\in\R^m\mid p_n(y)\geq0\}, \quad
M_- \coloneqq  \R^m\setminus M_+.
\]
Let $\chi_{+}, \chi_{-}$ denote the
characteristic functions of $M_{+}, M_{-}$ respectively.
Then define positive
Borel measures
\[
\dd\nu_{+}=\chi_{+}\dd\nu, \quad
\dd\nu_{-}=\chi_{-}\dd\nu, \quad
\dd\theta_{+} =   p_n\dd\nu_{+}, \quad
\dd\theta_{-} = - p_n\dd\nu_{-}.
\]
By definition, $\nu=\nu_++\nu_-$, so
\[
\int y_j^{2k}\,\dd\nu_+(y)\leq
\int y_j^{2k}\,\dd\nu(y)\quad \text{ for all }j,k\in\mathbb{N},
\]
whence $\nu_+$ satisfies the
Carleman condition~\reff{fix:suppnu=Q}.
Hence, so does $\theta_+$, again by
\cite[Corollary 14.22]{Scm}.
In particular, the measure $\theta_+$ is determinate.
Since $\dd\theta_+-\dd\theta_-=p_n\dd\nu$, we have
\[
\int y^\alpha \,\dd\theta_+(y)=
\int y^\alpha \,\dd\theta_-(y)
+ \int y^\alpha \,\dd\tau(y)
=
\int y^\alpha\, \dd(\theta_-+\tau)(y).
\]
Thus, by determinacy, $\theta_+=\theta_-+\tau$.
This yields
\[
0=\theta_+(M_-)\geq\theta_-(M_-)\geq0,
\]
so $\theta_-(M_-)=0$ and $\theta_-=0$.

Next, assume, for {the} sake of contradiction, that $p_n(y_0)<0$ for some $y_0\in\supp\nu$.
Then $-p_n(y)\geq\delta>0$ for all $y$
in a small ball $B$ around $y_0$.
This yields the contradiction\looseness=-1
\[
0=\theta_-(B)=\int (-p_n(y))\, \dd\nu_-(y)
=\int (-p_n(y))\, \dd\nu(y) \geq \delta \nu(B)>0,
\]
so $p_n\geq0$ on $\supp\nu$.
Finally, this again leads to the contradiction
\[
0>\ell_n(\phi)=\int \phi p_n\, {\rm d}\nu \geq0 ,
\]
which completes the proof.
\endproof

\subsection{The Positivstellensatz}

Now we consider the set
$K\subseteq\R^n$ as in \eqref{eq:K}.
Since $K$ is defined by the universal quantifier $y$ in $Q$,
one can {equivalently} write $K$ as the intersection\looseness=-1
\beq \label{K=intersect}
K \, =  \, \bigcap_{y\in Q}
\left \{x\in\R^n : g_1(x,y)\geq0,\ldots,g_s(x,y)\geq 0
\right \}.
\eeq
Clearly, $K$ is closed since each $g_i$ is a polynomial.
If the quantifier set $Q$ is semialgebraic,
then so is $K$ by Tarski's transfer principle \cite{BCR98}.
If $Q$ is not semialgebraic,
then $K$ may not be semialgebraic.

For notational convenience, denote
\[
g_0  \coloneqq   1, \quad g \coloneqq  (g_0, g_1, \ldots, g_s).
\]
For $f \in \rx$, if there exist SOS polynomials
$\tau_0, \tau_1, \ldots, \tau_s \in \sosxy$ such that
\beq\label{eq:fxQM}
f(\x) = \sum_{j=0}^s\int \tau_j(\x,y)g_j(\x,y)\, \dd \nu(y)
\eeq
then $f(x) \ge 0$ for all $x \in K$ {since $\supp{\nu} = Q$ by Assumption~\ref{ass:measure}}.
The set of all polynomials in $\R[\x]$ that can be represented as in \eqref{eq:fxQM} is
\beq\label{eq:QM}
\QM{g,\nu} \,  \coloneqq   \, \Bigg \{
\sum_{j=0}^s\int \tau_j(\x,y)g_j(\x,y)\, \dd \nu(y)  \mathrel{\Big|}
\text{each} \, \,  \tau_j  \in \sosxy
\Bigg \} .
\eeq
The set $\QM {g,\nu}$ is a convex cone in $\rx$.
It is called the {\it quadratic module} associated to $g$ and $\nu$, since
\[
\begin{array}{r}
	1\in \QM{g,\nu}, \quad
	\QM{g,\nu}+\QM{g,\nu}\subseteq \QM{g,\nu}, \\[1mm]
	\sosx \cdot \QM{g,\nu}\subseteq \QM{g,\nu}.
\end{array}
\]
Apparently, all polynomials in
$\QM{g,\nu}$ are nonnegative on $K$.
The Positivstellensatz concerns the reverse of this implication.
We start with the key Proposition \ref{prop:pos domain} stating that the positivity domain of %
$\QM{g,\nu}$ is $K$.

\begin{prop}\label{prop:pos domain}
    {Let $\nu$ be a  {nonnegative Borel} measure satisfying
   Assumption~\ref{ass:measure},}  then we have
\[   K=\{ x\in\R^n \mid f(x) \ge 0 \,\,  \forall f\in\QM{g,\nu} \}. \]
\end{prop}

\proof{Proof.}
By the definition \eqref{eq:QM}, every polynomial in $\QM{g,\nu}$
is nonnegative on $K$, whence
$K\subseteq\{ x\in\R^n \mid \forall f\in\QM{g,\nu}:\,f(x)\geq0\}=:\mathcal{D}$.

To establish the converse inclusion, assume
$\hat x\notin K$.	
Then there is a
$\hat y\in Q$ and a $j\in\{1,\ldots,s\}$ such that $g_j(\hat x,\hat y)<0.$
In a small open disk $B_{\eps_1}(\hat x,\hat y)$
of radius {$\eps_1>0$} about $(\hat x,\hat y)$
in $\R^{n+m}$, $g_j(x,y)  \leq -\la$ for some $\la>0.$
Consider a continuous
function $\phi$ positive on the open ball $B_{\frac{\eps_1}2}(\hat y)\subset\R^m$
and zero outside of {the closed ball} $\overline B_{\frac{\eps_1}2}(\hat y)$. Clearly,
\[
\psi(\x)\coloneqq\int_Q \phi(y) g_j(\x,y) \, {\rm d}\nu(y) \in \rx
\]
is negative at $\hat x$.

By Proposition \ref{prop:+approx}, there is a sequence $(\sigma_k)_k$ in $\sosy$ that converges
to $\phi$ in the $L^2$-norm.
Hence, for each $x$, as $k\to\infty$, we have
\[
\int_Q \sigma_k(y) g_j(x,y) \, {\rm d}\nu(y)
\longrightarrow  \int_Q \phi(y) g_j(x,y)\, {\rm d}\nu(y)=\psi(x).
\]
In particular, for $k$ large enough,
\[
f(\x) \coloneqq\int_Q \sigma_k(y) g_j(\x,y)\, {\rm d}\nu(y) \in\QM{g,\nu}
\]
is negative at $\hat x$. That is, $\hat x\not\in \mathcal{D}$,
whence $\mathcal{D}\subseteq K$ and we are done.
\endproof

\subsection{Bounded $K$}

In Positivstellens\"{a}tze,  we typically require that
$f>0$ on $K$ and the quadratic module associated to $K$ is archimedean.
Since $K$ is given by a universal quantifier over $y \in Q$,
we form the quadratic module $\QM{g,\nu}$ {and we} assume it is {\it archimedean}, i.e.,
there exists an integer $N>0$ such that
\[
N-\x_1^2 - \cdots - \x_n^2 \, \in \, \QM{g,\nu}.
\]
Clearly, the archimedeanness of $\QM{g,\nu}$ implies that $K$ is bounded
(so it is compact since it is closed).
Conversely, if $K$ is bounded,
we can generally assume $\QM{g,\nu}$ is archimedean,
because one can add the inequality $N-\sum\limits_{i=1}^{n} x_i^2\geq 0$ (no $y$)
to the description of the set $K$ as in \eqref{eq:K}.
The following is a generalization of the Putinar Positivstellensatz
to sets given by universal quantifiers.

\begin{theorem} \label{thm:univPsatz}
Let $K\subseteq\R^n$ be as in \eqref{eq:K} and
assume the measure $\nu$ satisfies Assumption~\ref{ass:measure}.
Suppose $\QM{g,\nu}$ is archimedean. For a polynomial $f\in\rx$,
if $f > 0$ on $K$, then we have $f \in \QM{g,\nu}$.
\end{theorem}

\proof{Proof.}
We shall apply \CRed{Jacobi's} \cite{Jacobi} strengthening of
Putinar's Positivstellensatz as presented in
\cite[Chapter 5]{Mar08}.
Consider the {archimedean} quadratic module $M=\QM{g,\nu}$. By Proposition \ref{prop:pos domain},
its positivity domain ($\mathcal K_M$ in Marshall's notation) is equal to $K$.
Hence, the \CRed{Jacobi-Putinar} Positivstellensatz
presented by Marshall in
\cite[Theorem 5.4.4]{Mar08} implies that every polynomial positive on $\mathcal K_M=K$
belongs to the quadratic module $M=\QM{g,\nu}$.
\endproof

Theorem~\ref{thm:univPsatz} clearly yields the following two corollaries.

\begin{corollary}\label{cor:univPsatz}
Let $K,\QM{g,\nu}$ be as in Theorem \ref{thm:univPsatz}
 {with $\nu$ satisfying Assumption~\ref{ass:measure}.}
Then the following are equivalent for $f\in\rx$:
\begin{enumerate}[label=(\roman*)]
	\item
	$f\geq0$ on $K$;
	\item
	for all $\eps>0$, $f+\eps\in \QM{{}g,\nu}$.
\end{enumerate}

\end{corollary}

\begin{corollary}\label{cor:empty}
Let $K,\QM{g,\nu}$ be as in Theorem \ref{thm:univPsatz}
 {with $\nu$ satisfying Assumption~\ref{ass:measure}.}
Then the following are equivalent:
\begin{enumerate}[label=(\roman*)]
\item
$K=\emptyset$;
\item
$-1\in \QM{{}g,\nu}$.
\end{enumerate}
\end{corollary}

\subsection{The non-archimedean case}

When the quadratic module $\QM{g,\nu}$ is not archimedean
(e.g., this is the case when $K$ is unbounded),
the conclusion of Theorem~\ref{thm:univPsatz}
may not hold. However,
Proposition \ref{prop:pos domain} allows us to get a
perturbation type Positivstellensatz
as in Lasserre-Netzer \cite{LN07},
for all (including unbounded) $K$.
For $r\in\N$, denote
\[
\Omega_r \coloneqq
\sum_{j=1}^n\sum_{k=0}^r \frac{\x_j^{2k}}{k!}\in\rx.
\]
We now have the following Positivstellensatz.

\begin{corollary}\label{cor:perturb}
Let $K\subseteq\R^n$ be as in \eqref{eq:K} and
 {with $\nu$ satisfying Assumption~\ref{ass:measure}.}
Then the following are equivalent for $f\in\rx$:\looseness=-1
\begin{enumerate}[label=(\roman*)]
	\item $f\geq0$ on $K$;
	
	\item
	for all $\eps>0$, there exists $r\in\N$ such that $f+\eps\Omega_r\in \QM{g,\nu}$.
\end{enumerate}
\end{corollary}
	
\proof{Proof.}
We shall apply a strengthening of the
Lasserre-Netzer perturbative Positivstellensatz {\cite[Corollary 3.6]{LN07}}
proved in \cite{KMV+} that can handle arbitrary constraints, and
is proved as a corollary of more general results on ``moment'' polynomials,
{i.e., polynomials in $\x$ and their formal moments with regard to a probability measure.}

Consider the constraint set $S =\QM{g,\nu}$.
In the notation of \cite{KMV+}, $K(S)=K$ and $Q(S)=\QM{g,\nu}$.
Now we simply apply \cite[Corollary 6.13]{KMV+}
(polynomials nonnegative on $K(S)$ are up to a perturbation as in (ii) contained in $Q(S)$)
to deduce
Corollary \ref{cor:perturb}.
\endproof

\subsection{Some illustrative examples}
\label{ssc:expo}

In the following examples, the measure $\nu$ is the classical Lebesgue measure.
Recall that $g_0 = 1$.

\begin{example}
Consider $f(\x)=- \x^3_1- \x^3_2 + \frac{1}{9} \x^2_1 \x_2 + \frac{1}{9} \x_1 \x^2_2  +
8 \x^2_1 + 8 \x^2_2$	and the set $K$ given as  {in \reff{K=intersect} with}
\begin{equation*}
	 {\begin{pmatrix}
	g_1(\x,\y)
		\\g_2(\x,\y)\\
	\end{pmatrix}} \, \coloneqq  \, \begin{pmatrix}
		1 - \x^2_1 \y^2_1 - \x^2_2 \y^2_2
		\\
		\x_1 \y^2_2 + \x_2 \y^2_1 -3\x_1\x_2 \y_1 \y_2\\
	\end{pmatrix}
	\]
	\[
	\mbox{and} \quad
	Q \, \coloneqq \, \left\{
	(y_1, y_2)
	\left \vert \begin{array}{c}
		y_1 + y_2 \leq 1, \\ y_1 \ge 0, \, y_2 \ge 0
	\end{array} \right.
	\right \}.
\end{equation*}
  {Note that the Lebesgue measure fulfills Assumption~\ref{ass:measure} when $Q$ is compact.}	
      A Positivstellensatz certificate for $f \in \QM{g,\nu}$ is
\begin{equation} \label{exm3.8:f:in:QM}
	f(\x)=\sum\limits_{i=0}^{2}\int_{0}^{1}\int_{0}^{1-y_2}\tau_i(\x,y)g_i(\x,y)\,\dd y_1\dd y_2,
\end{equation}
where the SOS polynomials $\tau_i(\x,\y)$ are
\begin{equation*}
	\begin{aligned}
		& \tau_0=(2 \x^2_1 - \x_1)^2 + (2 \x^2_2 - \x_2)^2 + 5 (\x_1 + \x_2)^2
		,\\			& \tau_1=60 (\x_1 \y_1 - \x_2 \y_2)^2,\\
		& \tau_2=20 (\x_1 \y_2 - \x_2 \y_1)^2 + 4(\x_1^2 + \x_2^2).\\
	\end{aligned}
\end{equation*}
       {One can check the representation \reff{exm3.8:f:in:QM}
      by a direct evaluation of integrals there.}
\end{example}

We remark that a Positivstellensatz certificate for $f \in \QM{g,\nu}$
can be computed numerically by solving a semidefinite program.
The following is such an example.

\begin{example}\label{Pos:exact:4}
Consider $f(\x)=\x^2_1\x_2-\x_1\x^2_2+\x^2_1+\x_2$
and the set $K$ given as  {in \reff{K=intersect} with}
\[
 {\begin{pmatrix}
		g_1(\x,\y)
		\\g_2(\x,\y)\\
\end{pmatrix}} \, \coloneqq \,\begin{pmatrix}
	\x^2_1\y_2+\x_2\y^2_1-\x_1+\x_2-\y_1\\
	\x_2\y^2_2-\x^2_2\y_1-\x_1\y_2+\x_2\y_1\\
\end{pmatrix}
\]
\[
  \mbox{and} \quad
Q = \{(y_1, y_2): |y_1|+|y_2| \leq 1\} .
\]
{Notice that the Lebesgue measure fulfills Assumption~\ref{ass:measure}
since $Q$ is compact. A  Positivstellensatz certificate $f \in \QM{g,\nu}$ is
\begin{equation} \label{exm3.9:finQM}
f(\x) = \sum\limits_{i=0}^{2}\int_{Q}\tau_i(\x,y)g_i(\x,y)\, \dd y,
\end{equation}
where $\tau_i(\x,\y)$ are SOS polynomials. We can represent them as
\[
\tau_0(\x) = [\x]_1^T X_0 [\x]_1, \quad \tau_1(\x,\y) = [\x,\y]_1^T X_1 [\x,\y]_1, \quad
\tau_2(\x,\y) = [\x,\y]_1^T X_2 [\x,\y]_1,
\]
where $X_0, X_1, X_2$ are symmetric positive semidefinite matrices.
By comparing coefficients of monomials of $\x$ in both sides of \reff{exm3.9:finQM},
we get a set of linear equations on $X_0, X_1, X_2$.
The matrices $X_0, X_1, X_2$ satisfying these conditions
can be found by solving a semi-definite program.
By using the software {\tt SeDuMi}, we obtained that}
\begin{equation*}
	\begin{aligned}
		& \tau_0=\begin{bmatrix}
			1\\\x_1\\\x_2
		\end{bmatrix}^{\top}
		\begin{bmatrix}
			0.0288 & 0.0988 & -0.0265 \\
			0.0988 & 0.3385 & -0.0909 \\
			-0.0265 & -0.0909 & 0.0244 \\
		\end{bmatrix}
		\begin{bmatrix}
			1\\\x_1\\\x_2
		\end{bmatrix}, \\
		& \tau_1= \begin{bmatrix}
			1\\\x_1\\\x_2\\\y_1\\\y_2\\
		\end{bmatrix}^{\top}
		\begin{bmatrix}
			0.0905 & -0.0988 & 0.0455 & 0.0865 & -0.2965 \\
			-0.0988 & 0.1080 & -0.0497 & -0.0945 & 0.3239 \\
			0.0455 & -0.0497 & 0.0229 & 0.0435 & -0.1492 \\
			0.0865 & -0.0945 & 0.0435 & 0.0827 & -0.2835 \\
			-0.2965 & 0.3239 & -0.1492 & -0.2835 & 0.9717
		\end{bmatrix}
		\begin{bmatrix}
			1\\\x_1\\\x_2\\\y_1\\\y_2\\
		\end{bmatrix},\\
		& \tau_2= \begin{bmatrix}
			1\\\x_1\\\x_2 \\ \y_1\\\y_2\\
		\end{bmatrix}^{\top}
		\begin{bmatrix}
			0.3787 & 0.4577 & 0.0895 & 0.5813 & -0.1114 \\
			0.4577 & 1.9459 & -0.0505 & 1.0328 & -0.1879 \\
			0.0895 & -0.0505 & 0.0392 & 0.0998 & -0.0203 \\
			0.5813 & 1.0328 & 0.0998 & 0.9704 & -0.1836 \\
			-0.1114 & -0.1879 & -0.0203 & -0.1836 & 0.0348
		\end{bmatrix}
		\begin{bmatrix}
			1\\\x_1\\\x_2\\ \y_1\\\y_2\\
		\end{bmatrix}.
	\end{aligned}
\end{equation*}
The above matrices in the middle are all positive semidefinite.
For neatness, only four decimal digits are shown
(the errors for matching coefficients are in the order of $10^{-11}$).\looseness=-1
\end{example}

\begin{example}  \label{exm3.10:Q=Z}
Consider $f(\x)= 4 - \x_1^2 - \x_2^2 $ and the set $K$ given as in \reff{K=intersect} with
\[
 g_1(\x,\y)   \, \coloneqq \, 
	 4\y^4 -1 - \y^2(2\y^2-1) \x_1^2 - \y^2(2\y^2+1) \x_2^2 
\]
\[
\mbox{and} \quad
Q = \Z_+ =  \{ 1, 2, \ldots \} .
\]
We select the measure $\nu$ supported on $Q$ such that
\[
  \nu( \{ k \} ) \, = \, \frac1e\cdot\frac1{k!}, \quad \text{for} \quad k= 1, 2, \ldots .
\]
One can directly calculate that
\[
 \int y^j \mt{d} \nu(y) =\frac1e \sum_{k=1}^{\infty} \frac{  k^j}{k!} = \left\{\begin{array}{ll}
 	1-\frac{1}{e} & \quad\text{if}\quad j=0,\\
 	B_j & \quad\text{if}\quad j\geq 1.
 \end{array}\right.
\]
In the above, $B_j$ denotes the $j$th Bell number (\cite[Section II.3]{FS09} or \cite[p.~82]{stanley}),  
which counts the number of partitions of the set $[j]=\{1,\ldots,j\}$.
It is interesting to remark that $B_j \leq j!$, which can be seen as follows.
We assign to each partition $[j]=S_1\cup\cdots\cup S_k$ a different permutation as follows:
sort each $S_i$ in increasing order, and relabel $S_i$ so that the lowest number in
$S_i$ is smaller than the lowest number in $S_{i+1}$.
Then each $S_i$ yields a permutation when it is viewed as a cycle,
and the product of the disjoint cycles assigned to the $S_i$'s is a permutation of $[j]$.
We now claim that $\nu$ satisfies the Carleman condition \eqref{fix:suppnu=Q}. 
Indeed, by Stirling's approximation for factorials
(see, e.g., \cite[Section I.2]{FS09} or \cite{Rob55}),
$n!\sim \sqrt{2\pi n}\left(\frac ne\right)^n $,
we have
\[
 \sum_{d=1}^{\infty} \Big( \int y^{2d} \mt{d} \nu(y) \Big)^{ -\frac{1}{2d} } \ge
 \sum_{d=1}^{\infty}  \big( (2d)! \big)^{-\frac{1}{2d} }  \sim
  \sum_{d=1}^{\infty}\frac{e}{2d}  \big( \sqrt{4 \pi d}  \big)^{-\frac{1}{2d} }
  = \infty  .
\]
The set $K$ is the intersection of infinitely many ellipses.
The Positivstellensatz certificate $f \in \QM{g,\nu}$ is
\begin{equation*}
	f(\x) = \int_{Q}\tau_0(\x,\y)  \dd \y + \int_{Q}\tau_1(\x,\y)g_1(\x,\y)\, \dd \y
\end{equation*}
where $\tau_0(\x,\y), \tau_1(\x,\y)$ are SOS polynomials.
By solving a semidefinite program, we can get
\begin{equation*}
	\begin{aligned}
		& \tau_0 \approx
             0.1671+0.6336\y+0.0990\x_2^2+0.6005\y^2, \\
				& \tau_1  \approx
         0.0053-0.0100\y+0.0047 \y^2 	.
	\end{aligned}
\end{equation*}
For neatness, only four decimal digits are shown in the above
(the errors for matching coefficients are in the order of $10^{-12}$).
\end{example}

In Example~\ref{exm3.10:Q=Z}  the set $K$ is compact.
Now consider  the same $Q = \Z_+$ and the measure $\nu$ as in the above example.  If
$g(\x,\y) = \x_2 - 2 \y \x_1  + \y^2$, then the set
$K\supseteq\{x_1 \le 0, x_2 \ge 0 \} \cup \{x_1 \ge 0, x_2 - x_1^2 \ge 0\}$
is an unbounded convex region, and
the quadratic module $\QM{g, \nu}$ cannot be archimedean.

\section{Moment problems with universal quantifiers}
\label{sc:Mom}

This section considers $K$-moment problems for
the set $K$ given by a universal quantifier as in \reff{eq:K}.
Recall from the introduction that
a multisequence $z = (z_\al)_{\al \in \N^n}$
gives rise to the {Riesz functional}
\begin{equation} \label{eq:rieszFun}
	\mathscr{R}_z:\rx\to\R, \quad \x^\al \mapsto z_\al .
\end{equation}
Equivalently,
\[
\mathscr{R}_z \Big( \sum_\af f_\af \x^\af \Big) = \sum_\af f_\af z_\af,
\]
where $f_\af\in\R$ are {real} coefficients.
If $\mu$ is a representing measure for $z$, then
\[
\mathscr{R}_z(f) = \int f(x)\, \dd \mu(x)
\quad  \text{for all} \quad  f \in \rx.
\]
If  $\supp{\mu} \subseteq K$,
then $z$ must satisfy
\begin{equation} \label{Rz(f)>=0:allindx}
	\mathscr{R}_z(f) \, \ge \, 0
	\quad  \text{for all} \quad f \in \rx: \, f|_K \ge 0.
\end{equation}

The multisequence $z$ is said to be $K$-positive if \reff{Rz(f)>=0:allindx} holds.
Clearly, being $K$-positive is a necessary condition
for $z$ to have a $K$-representing measure.
When the set $K$ is closed (this is the case if $K$ is given as in \reff{eq:K}),
being $K$-positive as in \reff{Rz(f)>=0:allindx} is also sufficient
for $z$ to be a $K$-moment sequence.
This is a classical result of Riesz and Haviland.
The reader is referred to the surveys
\cite{Ber,Fia16} and books \cite{Akh,Scm}
for more details about classical moment problems.

Using the quadratic module
$\QM{g,\nu}$ introduced in  \eqref{eq:QM},
we have the following characterization of a $K$-moment sequence.

\begin{theorem}\label{theorem:mp}
	Let $K\subseteq\R^n$ be as in \eqref{eq:K} and suppose
	the measure $\nu$ satisfies Assumption~\ref{ass:measure}.
	Assume the quadratic module $\QM{g,\nu}$ is archimedean.
	Then the multisequence $z$ is a $K$-moment sequence
	if and only if the Riesz functional
	$\mathscr{R}_z \ge 0$ on $\QM{g,\nu}$, i.e.,
	$\mathscr{R}_z(f) \ge 0$ for all $f \in \QM{g,\nu}$.
\end{theorem}
  \proof{Proof.}
\noindent $(\Rightarrow)$ If $f \in \QM{g,\nu}$,
then $f|_K \ge 0$. Hence \reff{Rz(f)>=0:allindx}
implies $\mathscr{R}_z \ge 0$ on $\QM{g,\nu}$.
Conversely, if $\mathscr{R}_z \ge 0$ on $\QM{g,\nu}$,
then $\mathscr R_z$ is also nonnegative on each
$f\in\R\cx$ that is nonnegative on $K$ by Theorem \ref{thm:univPsatz} or Corollary \ref{cor:univPsatz}.
The implication $(\Leftarrow)$ now follows by the Riesz-Haviland theorem mentioned above.
\endproof

As pointed out by one of the referees, an alternate proof of Theorem \ref{theorem:mp} can be given using the Jacobi-Putinar Positivstellensatz, cf.~\cite[Section 5.6]{Mar08}.

\begin{corollary}\label{cor:mp1}
Let $K\subseteq\R^n$ be as in \eqref{eq:K} and suppose
the measure $\nu$ satisfies Assumption~\ref{ass:measure}.
Then the multisequence $z$ is a $K$-moment sequence
if and only if the Riesz functional
$\mathscr{R}_z$ satisfies
$\mathscr R_z(f)\geq 0$ for all $f \in \rx$
 {with the following property: for all $\eps >0$ there exists $r\in\N$ with
$f+\eps \Omega_r\in\QM{g,\nu}$.}
\end{corollary}

\proof{Proof.}
 {\noindent $(\Rightarrow)$ is obvious since a polynomial $f$ satisfying the perturbation condition
in the statement of the corollary, is nonnegative on $K$ (cf.~Proposition \ref{prop:pos domain}
or Corollary \ref{cor:perturb}). For the converse implication $(\Leftarrow)$
note that $f$ is nonnegative on $K$
if and only if it satisfies this perturbation condition (again by Corollary \ref{cor:perturb}).
The conclusion now follows by the Riesz-Haviland theorem.}
\endproof

\begin{remark}\label{rem:ref2}
As pointed out by one of the referees, a strengthening of Theorem \ref{theorem:mp} holds.
Namely, assume $K\subseteq\R^n$ is as in \eqref{eq:K}, and
the measure $\nu$ satisfies Assumption~\ref{ass:measure}.
Then, if the multisequence $z$ itself satisfies the multivariate Carleman condition, 
then it is a $K$-moment sequence iff $\mathscr{R}_z \ge 0$ on $\QM{g,\nu}$. 
Indeed, the proof is essentially the same as that of Theorem \ref{theorem:mp}, 
but at the final step one applies \cite[Theorem 3.16]{IKKM22}, 
a far reaching extension of the Nussbaum theorem.
\end{remark}
 
In the sequel, we determine concrete conditions on $z$ for
$\mathscr{R}_z \ge 0$ on $\QM{g,\nu}$.

\subsection{Moment and localizing matrices}
\newcommand{\mvec}[1]{\mathbf{#1}}

The multisequence $z = (z_{\al})_{\al \in \N^n }$ gives rise to the infinite matrix
\[
\mathrm{H}[z] \, \coloneqq \, (z_{\al + \be})_{\al,\be \in \N^n }.
\]
That is, $\mathrm{H}[z]$ is the matrix labelled by
nonnegative integer vectors $\al, \be \in \N^n$ and
\[
\mathrm{H}[z]_{\al,\be} \,= \, z_{\al+\be}
\]
for all $\al,\be$.
It is called the {\it moment matrix} or {\it multivariate Hankel matrix}
of the multisequence $z$.
For a vector $\mvec{u}=(u_\al)_{\al \in \N^n}$ with  finitely many nonzero entries,
we have\looseness=-1
\[
\mvec u^\T \mathrm{H}[z] \mvec u = \mathscr{R}_z\Big( u(\x)^2 \Big),
\quad \text{where} \quad
u(\x) = \sum_\al u_\al \x^\al.
\]
Hence, if $\mathscr{R}_z \ge 0$ on $\Sigma^2[\x]$, then $\mathrm{H}[z] \succeq 0$.
For a degree $k$, we denote the truncation
\[
\mathrm{H}^{(k)}[z]\, \coloneqq \, (z_{\al + \be})_{\al,\be \in \N_{k}^n }.
\]
One can easily verify that $\mathrm{H}^{(k)}[z] \succeq 0$
if $\mathscr{R}_z \ge 0$ on $\sosx\cap\R\cx_{2k}$.

Next we give localizing matrices for the quadratic module $\QM{g,\nu}$.
For a given multisequence $z$,
$\mathscr{R}_z \Big( \int p(\x,y)^2 g_j(\x,y)\, \dd \nu(y)\Big)$
is a quadratic form in the vector of coefficients of $p(\x,\y)$.
For convenience, we use $\mvec p$ to denote the vector of coefficients of $p(\x,\y)$.
Let $\mathrm{L}_{\nu,g_j}^{(k,l)}[z]$ be the matrix
associated to this quadratic form. Here superscripts $k,l$
denote degree bounds on $\x$ and $\y$, respectively, so that
\[
\mathscr{R}_z \Big( \int p(\x,y)^2 g_j(\x,y)\, \dd \nu(y) \Big)
\, = \,  \mvec p^\T \Big( \mathrm{L}_{\nu, g_j}^{(k,l)}[z] \Big) \mvec p,
\]
for all $p(\x,\y) \in \rxy$ with degrees
\beq \label{deg:p2gj<=2k2l}
\deg_\x( p(\x,\y)^2 g_j(\x,\y) ) \le 2k, \quad
\deg_\y( p(\x,\y)^2 g_j(\x,\y) ) \le 2l.
\eeq
Explicit expressions for $\mathrm{L}_{\nu,g_j}^{(k,l)}[z]$ can be given as follows.
For convenience, denote %
\beq \label{kprm:lprm}
k'  \coloneqq   k - \lceil \deg_\x(  g_j(\x,\y) )/2 \rceil, \quad
l'  \coloneqq   l - \lceil \deg_\y(  g_j(\x,\y) )/2 \rceil .
\eeq
Then we can write
\[
p(\x,\y) = \mvec p^\T([\x]_{k'} \otimes [\y]_{l'})
\]
where $[\x]_k$ denotes the vector of all monomials in $\x$
of degrees at most $k$, and likewise for $[\y]_l$.
The constraining polynomial $g_j(\x,\y)$ can be written in the form
\[
g_j(\x,\y) = \sum_i g_{ji}(\x) h_{ji}(\y),
\]
for some polynomials $g_{ji} \in \rx$ and $h_{ji} \in \ry$.
Then, one can see that
\[
\begin{array}{rl}
	& \mathscr{R}_z \Big( \int p(\x,y)^2 g_j(\x,y)\, \dd \nu(y) \Big) \\
	=& \mvec p^\T \Big( \mathscr{R}_z  \int   g_j(\x,y)
	([\x]_{k'} \otimes [y]_{l'}) ([\x]_{k'} \otimes [y]_{l'} )^\T\, \dd \nu(y) \Big) \mvec p \\
	=&\mvec p^\T \Big( \mathscr{R}_z  \int   g_j(\x,y)
	([\x]_{k'} [\x]_{k'}^\T) \otimes [y]_{l'} [y]_{l'}^\T\, \dd \nu(y) \Big) \mvec p \\
	= & \mvec p^\T \Big( \sum\limits_i
	\big( \int  h_{ji}(y) [y]_{l'} [y]_{l'}^\T \dd\, \nu(y) \big) \otimes
	\mathscr{R}_z \big(  g_{ji}(\x) [\x]_{k'} [\x]_{k'}^\T \big) \Big)  \mvec p .
\end{array}
\]

(In the above, when $\mathscr{R}_z$ is applied to a matrix,
it means that it is applied entrywise, for convenience of notation.)
Denote the matrices
\beq \label{mat:YL}
Y_{\nu,h_{ji}}^{(l')} \,  \coloneqq   \, \int  h_{ji}(y) [y]_{l'} [y]_{l'}^\T \, \dd \nu(y), \quad
\mathrm{L}_{g_{ji}}^{(k')}[z] \,  \coloneqq   \, \mathscr{R}_z
\Big(  g_{ji}(\x) [\x]_{k'} [\x]_{k'}^\T  \Big) .
\eeq
Then we get the expression
\beq \label{express:Lhji[z]}
\mathrm{L}_{\nu,g_j}^{(k,l)}[z] \, \coloneqq \,
\sum_i Y_{\nu,h_{ji}}^{(l')} \otimes \mathrm{L}_{g_{ji}}^{(k')}[z] .
\eeq
Note that $k',l'$ are the degrees defined in \reff{kprm:lprm}.
Observe that $\mathrm{L}_{g_{ji}}^{(k')}[z]$ is the localizing matrix
for the polynomial $g_{ji} \in \rx$, and is independent of $\nu$.
Similarly, the matrices $Y_{\nu,h_{ji}}^{(l')}$ are independent of $z$.
In particular, for $g_0 = 1$, we get
\beq
\mathrm{L}_{\nu, 1}^{(k,l)}[z]  \,= \,
\Big( \int   [y]_{l} [y]_{l}^\T \dd \nu(y) \Big)
\otimes \mathrm{H}^{(k)}[z].
\eeq

\subsection{The full moment problem}

We give a full characterization for $K$-moment sequences
when $K$ is defined by universal quantifiers.

\begin{theorem}\label{thm:shift}
Let $K\subseteq\R^n$ be as in \eqref{eq:K} and assume
the measure $\nu$ satisfies Assumption~\ref{ass:measure}.
Then, for a multisequence $z$,
we have $\mathscr{R}_z \ge 0$ on $\QM{g,\nu}$
if and only if for all $j = 0, 1, \ldots, s$,
\beq \label{Lgj[z]>=0:allk}
\mathrm{L}_{\nu, g_j}^{(k,l)}[z] \, \succeq \, 0,
\quad
k=1,2,\ldots, \;
l=1,2,\ldots .
\eeq
Moreover, when $\QM{g,\nu}$ is archimedean,
then $z$ is a $K$-moment sequence
if and only if it satisfies \reff{Lgj[z]>=0:allk}.
\end{theorem}
\proof{Proof.}
Observe that $\mathscr{R}_z \ge 0$ on $\QM{g,\nu}$ if and only if
\[
\mathscr{R}_z \Big( \int p(\x,y)^2 g_j(\x,y)\, \dd \nu(y) \Big)  \ge 0
\]
for all $j$ and for all $p(\x,\y) \in \rxy$.
When $p$ is restricted to have degrees as in \reff{deg:p2gj<=2k2l},
then \reff{Lgj[z]>=0:allk}
follows from the definition of
$\mathrm{L}_{\nu,g_j}^{(k,l)}[z]$
for all $k$ and $l$.  When $\QM{g,\nu}$ is archimedean,
the last statement follows from  Theorem~\ref{theorem:mp}.\hspace{1cm}
{\endproof}

When $K$ is given  without quantifiers,
there is a classical flat extension theorem \cite{CF96,CurFia05}
that recognizes $K$-moment sequences.
Here, we give a similar flat extension theorem for
sets $K$ defined with universal quantifiers. Let
\beq \label{deg:dg}
d_g \, \coloneqq \, \max\{1, \deg_{\x}(g)  \}.
\eeq

\begin{theorem} \label{thm:FE:full}
Let $K\subseteq\R^n$ be as in \eqref{eq:K} and assume
the measure $\nu$ satisfies Assumption~\ref{ass:measure}.
Let $z$ be a multisequence satisfying \reff{Lgj[z]>=0:allk}.
If there exists $k \ge d_g$ such that
\[
r \coloneqq
\rank  \mathrm{H}^{(k-d_g)}[z] \, = \, \rank \mathrm{H}^{(k)}[z],
\]
then $z$ admits an $r$-atomic measure $\mu$ supported in $K$
and  $\mu$ is the unique representing measure for $z$.
\end{theorem}
\proof{Proof.}
By the flat extension theorem \cite{CF96,CurFia05}, we know that $z$
admits an $r$-atomic representing measure, say, $\mu$.
Moreover, the $\mu$ is the unique representing measure for $z$.
Since $z$ satisfies \reff{Lgj[z]>=0:allk},
$\mathscr{R}_z \ge 0$ on $\QM{g,\nu}$.
Pick an arbitrary $f \in \QM{g,\nu}$,
then \reff{Lgj[z]>=0:allk} implies that
\[
\mathrm{L}_{f}^{(k)}[z] \, \succeq \, 0
\]
for all $k=1,2, \ldots.$
As in \cite{CF96,CurFia05}, we have $\supp{\mu} \subseteq \{x : f(x) \ge 0 \}$.
As this holds for all $f \in \QM{g,\nu}$, Proposition
\ref{prop:pos domain}
implies that $\supp{\mu} \subseteq K$.
\endproof

\subsection{The truncated moment problem}

Now we consider  $w = (w_\al)_{\al \in\N_{2d}^n}$, a truncated multisequence
of even degree $2d$. We look for concrete conditions under which
$w$ is a $K$-moment sequence, with a representing measure $\mu$ supported in $ K$.
As in the above calculations, for $w$ to be a $K$-moment sequence,
it must satisfy
\beq \label{Rw>=0}
\mathscr{R}_w \Big( \int p(\x,y)^2 g_j(\x,y)\, \dd \nu(y) \Big)  \ge 0
\eeq
for all $j$ and for all $p(\x,\y) \in \rxy$ with the degree
\[
\deg_{\x}  \big( p(\x,\y)^2 g_j(\x,\y) \big) \le 2d.
\]
Note that \reff{Rw>=0} is equivalent to
\beq \label{Lgj[z]>=0:k=1:d}
\mathrm{L}_{\nu, g_j}^{(k,l)}[w] \, \succeq \, 0, \quad
k=1,2,\ldots d, \; \;
l=1,2,\ldots .
\eeq

The following is a generalization of
the flat extension theorem in \cite{CF96,CurFia05}.

\begin{theorem} \label{thm:truncMP}
Let $K\subseteq\R^n$ be as in \eqref{eq:K} and assume
the measure $\nu$ satisfies Assumption~\ref{ass:measure}.
Let {$w \in \re^{\N^n_{2d} }$} be a truncated multisequence satisfying \reff{Lgj[z]>=0:k=1:d}.
If there exists a positive integer {$k \le d-d_g$} such that
\beq \label{FEC:w}
r  \, \coloneqq \,  \rank  \mathrm{H}^{ {(k)} }[w] \, = \,
\rank \mathrm{H}^{ {(d)} }[w],
\eeq
then $w$ admits an $r$-atomic measure $\mu$ supported in $K$
and $\mu$ is the unique representing measure for $w$.
\end{theorem}
\proof{Proof.}
By the flatness condition \reff{FEC:w}, the truncated multisequence $w$
can be extended to a full multisequence
$z = (z_\al)_{\al \in \N^n}$ with $\rank \mathrm  H[z]=r$
that represents an $r$-atomic measure $\mu$
       {(see Corollary~1.4 of \cite{Lau05}).
Moreover, $\mu$ is the unique representing measure for $w$ (and also for $z$).
This can be implied by Theorems~1.2 and 1.6 of \cite{Lau05} (also see \cite{CF96}).}
It now remains to show that $\supp{\mu} \subseteq K$.
For all $a(\x) \in \rx_d$ and $b(\y) \in \ry_l$, it holds that
\[
\mathscr{R}_w \Big( a(\x)^2  \int b(y)^2 g_j(\x,y)\, \dd \nu(y) \Big)
=
\mathscr{R}_w \Big( \int a(\x)^2 b(y)^2 g_j(\x,y)\, \dd \nu(y) \Big)  \ge 0.
\]
Since $z$ is an extension of $w$, {which is represented by $\mu$}, we get
\[
\mathscr{R}_z \Big( f(\x)^2  \int b(y)^2 g_j(\x,y)\, \dd \nu(y) \Big)
\ge 0
\]
for all $f \in \rx$. This implies that for each $j$,
\[
\supp{\mu} \subseteq \left\{ x\in\R^n \Big \vert  \int b(y)^2 g_j(x,y)\, \dd \nu(y)
\ge 0 \right\}.
\]
The above is true for all $b \in \ry_l$ and $l=1, 2, \ldots$.
Hence,  as in the proof of Proposition \ref{prop:pos domain}, we can infer
that the intersection over $j$ of the
right-hand side sets in the above equation is equal to $K$. That is,
$\supp{\mu} \subseteq K$, which completes the proof.\looseness=-1
\endproof

We remark that the rank condition \reff{FEC:w} implies that
the truncated multisequence $w$ admits a unique $r$-atomic representing
measure $\mu$, say, $w = \lmd_1 [u_1]_{2d} + \cdots + \lmd_r [u_r]_{2d}$,
for distinct points $u_1, \ldots, u_r$ and positive scalars
$\lmd_1, \ldots, \lmd_r$, as in \cite{CF96,CurFia05}.
The condition \reff{Lgj[z]>=0:k=1:d}
ensures that all $u_1, \ldots, u_r \in K$.
Note that \reff{Lgj[z]>=0:k=1:d} requires it to hold
for all $l = 1, 2, \ldots$. If this is not checkable,
one can verify $u_i \in K$ by
checking nonnegativity of $g(u_i, \y)$ on $Q$.
The following is such an example.

\begin{example}\label{tms_1}
For the  set
\[
K =\left\{x\in\R^2\mid 1-x^Ty\geq 0\quad \forall y\in\mathbb{R}^2: y^4_1+y^4_2\leq 1\right\},
\]
we consider the truncated multisequence $w \in \re^{ \N_4^2 }$ given such that
\[	
\mathrm{H}^{(2)}[w] \, = \,
\begin{bmatrix*}[r]
	3 & 0 & \frac{2}{3} & \frac{2}{3} & -\frac{5}{18} & \frac{17}{18} \\[6pt]
	0 & \frac{2}{3} & -\frac{5}{18} & -\frac{2}{9} & \frac{13}{54} & \frac{7}{108} \\[6pt]
	\frac{2}{3} & -\frac{5}{18} & \frac{17}{18} & \frac{13}{54} & \frac{7}{108} & \frac{8}{27} \\[6pt]
	\frac{2}{3} & -\frac{2}{9} & \frac{13}{54} & \frac{2}{9} & -\frac{23}{162} & \frac{61}{324} \\[6pt]
	-\frac{5}{18} & \frac{13}{54} & \frac{7}{108} & -\frac{23}{162} & \frac{61}{324} & -\frac{17}{648} \\[6pt]
	\frac{17}{18} & \frac{7}{108} & \frac{8}{27} & \frac{61}{324} & -\frac{17}{648} & \frac{209}{648}
\end{bmatrix*}.
\]
One can check that $\rank \mathrm{H}^{(1)}[w]=\rank \mathrm{H}^{(2)}[w]=3$,
so the condition \reff{FEC:w}
of Theorem \ref{thm:truncMP}
holds.  As in \cite{CurFia05},
we obtain $w=[u_1]_4+[u_2]_4+[u_3]_4$ for the points
\[
u_1=\left(-\frac{2}{3},\frac{1}{2}\right), \quad
u_2=\left(\frac{1}{3},\frac{2}{3}\right),  \quad
u_3=\left(\frac{1}{3},-\frac{1}{2}\right) .
\]
It is easily seen (e.g., by H\"older's inequality)
that these three points belong to the set $K$.
\end{example}

\section{Semi-Infinite Optimization}
\label{sc:SIP}

An important application of Positivstellens\"{a}tze
and moment problems with universal quantifiers
is to solve semi-infinite optimization.
Consider the semi-infinite program (SIP):
\begin{equation} \label{primal-SIP}
	\left\{
	\begin{array}{lll}
		\min\limits_{x\in X}& f(x)\\
		\st    &  g(x,y)\geq 0\quad\forall y\in Q.  \\
	\end{array}
	\right.
\end{equation}
The constraining function $g$ in \reff{primal-SIP}
is the $s$-dimensional vector of polynomials,
\[
g(\x,\y) \, \coloneqq \, \big( g_1(\x,\y), \ldots,  g_s(\x,\y) \big),
\]
$f\in\R\cx$, and
$X\subseteq\mathbb{R}^n$
is another given constraining set
that does not depend on $y\in\R^m$.
We assume $X$ is given as
\begin{equation}  \label{set:X}
	X=\left\{ x\in\mathbb{R}^n \;\middle\vert\;
	c_i(x)= 0 \,(i\in \mathcal{I}),\;
	c_j(x)\geq 0 \,(j\in\mathcal{J})
	\right\}  .
\end{equation}
Here, all $c_i, c_j$ are polynomials in $\x$ and
$\mathcal{I},\,\mathcal{J}$
are disjoint finite label sets. For convenience of notation,
we denote the polynomial tuples:
\[
c_{eq} = ( c_i)_{ i\in \mathcal{I} }, \quad
c_{in} = ( c_j)_{ j\in \mathcal{J} }.
\]

Semi-infinite optimization has broad applications,
such as Chebyshev approximation \cite{lopez2007semi} and
robustness support vector machines \cite{xu2009robustness}.
Classical methods for solving semi-infinite optimization
include Karush–Kuhn–Tucker multipliers \cite{stein2012adaptive},
discretization methods \cite{djelassi2017hybrid},
and Moment-SOS relaxations \cite{HuNie23,WangGuo14}.
In this section, we show how to use Positivstellens\"{a}tze
and moment problems with universal quantifiers
to solve SIPs.\looseness=-1

As before, we let $\nu$ be a nonnegative Borel measure on $\re^m$
satisfying Assumption~\ref{ass:measure}.
We assume the moments $\int_Q y^\alpha\, \dd \nu(y)$ are available.
Then truncations for given degrees of the quadratic module $\QM{g,\nu}$
can be represented by semidefinite programs.

\begin{prop}\label{prop:dual}
Let $K$ be as in \eqref{eq:K} and let $\nu$ be a
Borel measure satisfying Assumption \ref{ass:measure}.
Assume %
that the quadratic module $\qmod{g,\nu} + \ideal{c_{eq}} + \qmod{c_{in}}$ is archimedean.
Then the optimal value $f_{\min}$ of \reff{primal-SIP}
is equal to the optimal value of the following optimization problem
\beq \label{min:Rz(f):onC}
\left\{ \begin{array}{rl}
	\min & \mathscr{R}_z(f) \\
	\st  & \mathscr{R}_z \geq 0 \quad \mbox{on} \quad
	\qmod{g,\nu} + \ideal{c_{eq}} + \qmod{c_{in}},  \\
	& \mathscr{R}_z(1) = 1, \quad z \in \re^{ \N^n }.
\end{array} \right.
\eeq
\end{prop}

\proof{Proof.}
The feasible set of \reff{primal-SIP} is $X \cap K$,
where $K$ is as in \eqref{eq:K}. Note that
\[
X \cap K =\left\{ x\in\mathbb{R}^n \;\middle\vert\;
\begin{array}{c}
	\begin{bmatrix} c_{eq}(x) \\  -c_{eq}(x) \\ c_{in}(x) \\  g(x, y) \end{bmatrix}
	\ge 0 \; \forall y \in Q
\end{array}  \right\}.
\]
The polynomials $c_i$ can also be viewed as depending on $\y$ trivially.
Observe that
\[
\qmod{(c_{eq}, -c_{eq}, c_{in}, g), \, \nu}  \, = \,
\qmod{g,\nu} + \ideal{c_{eq}} + \qmod{c_{in}} .
\]
Since $\qmod{g,\nu} + \ideal{c_{eq}} + \qmod{c_{in}}$ is archimedean,
the set $X \cap K$ is bounded (cf.~Proposition \ref{prop:pos domain}).
Thus the optimal value $f_{\min}$ is finite, i.e.,  $f_{\min} \in \re$.
Hence, $f_{\min}$ equals the minimum value of
the expectation $\int_{X\cap K} f(x)\, \dd \mu(x)$, over all probability measures
$\mu$ supported in $X \cap K$.
When $z$ is a multisequence satisfying the constraints in \reff{min:Rz(f):onC},
Theorem~\ref{theorem:mp} implies that
$z$ is the moment sequence of such a probability measure $\mu$.
Therefore, $f_{\min}$ is also the minimum value of \reff{min:Rz(f):onC}.
\endproof

Proposition \ref{prop:dual} can be used to give
Moment-SOS type relaxations for solving
the semi-infinite optimization \reff{primal-SIP}.
The full multisequence $z \in \re^{ \N^n }$ can be
approximated by its truncations
\[ w = (z_\af)_{ \al \in \N_{2k}^n } , \]
for a degree $k$. Note that $\mathscr{R}_w \ge 0$ on $\QM{g,\nu}_{2k}$
if and only if
\[ \mathrm{L}_{\nu, g_j}^{(k,l)}[w] \succeq 0 \]
for all $l = 1, 2, \ldots$ (cf.~Theorem \ref{thm:shift}).
The constraining polynomials $c_j$ do not depend on $\y$, so
\[
\mathrm{L}_{\nu, c_j}^{(k,l)}[w]  \, = \,
\Big( \int 1\, \dd \nu(y) \Big) \cdot \mathrm{L}_{c_j}^{(k)}[w] .
\]
In computational practice, we typically scale $\nu$
so that $\int 1\, \dd \nu(y) = 1$, whence
$\mathrm{L}_{\nu, c_j}^{(k,l)}[w]  = \mathrm{L}_{c_j}^{(k)}[w]$.
It is also interesting to note that
$\mathscr{R}_z \geq 0$ on
$\qmod{g,\nu} + \ideal{c_{eq}} + \qmod{c_{in}}$
if and only if
$\mathscr{R}_z \geq 0$ on
each of the
$\qmod{g,\nu}, \, \ideal{c_{eq}}, \, \qmod{c_{in}}$. Moreover, $\mathscr{R}_z \geq 0$ on
$\ideal{c_{eq}}$ if and only if
$\mathscr{R}_z \equiv 0$ on
$\ideal{c_{eq}}$, since $\ideal{c_{eq}}$ is a subspace of $\rx$.
Note that $\mathscr{R}_z \equiv 0$ on $\ideal{c_{eq}}_{2k}$
is equivalent to
$\mathrm{L}_{c_i}^{(k)}[w] = 0$, for each $i \in \mathcal{I}$.\looseness=-1

Suppose $\deg(c_i) \le 2k$ for each $i$.
Let $\mathscr{V}_{c_i}^{(2k)}[w]$
denote the vector such that
\beq \label{locvec:Vci}
\mathscr{R}_w ( c_i(\x) u(\x) ) \, = \,
\big( \mathscr{V}_{c_i}^{(2k)}[w] \big)^T \mvec u
\eeq
for all $u(\x) \in \R\cx_{2k - \deg(c_i)}$.
The $\mathscr{V}_{c_i}^{(2k)}[w]$  is called the
{\it localizing vector} of the polynomial $c_i$,
generated by the truncated multisequence $w$.
It is important to observe that
$\mathscr{V}_{c_i}^{(2k)}[w] = 0$
if $w$ has a representing measure supported on $c_i(x) = 0$.

To get a finite dimensional optimization problem,
we choose a finite value for $l$, e.g., $l= k$. Recall that
\[
\mathscr{R}_w(f) = \sum_\af f_\af w_\af
\quad \text{for} \quad  f(\x) = \sum_\af f_\af \x^\af .
\]
In particular, $w_0 = \mathscr{R}_w(1) = 1$.
Therefore, the $k$th order truncation of \reff{min:Rz(f):onC} is
\beq \label{min:Rz(f):order=k}
\left\{ \begin{array}{rl}
	\gamma_k  \coloneqq   \min & \sum\limits_\af f_\af w_\af \\[1mm]
	\st  & \mathscr{V}_{c_i}^{(2k)}[w] = 0  \,(i \in \mathcal{I}), \\[1mm]
	& \mathrm{L}_{c_j}^{(k)}[w] \succeq 0  \,(j \in \mathcal{J}),  \\[1mm]
	& \mathrm{L}_{\nu, g_j}^{(k,k)}[w] \succeq 0  \,(j=0,1,\ldots, s),  \\ [1mm]
	& w_0 = 1, \quad w \in \re^{ \N_{2k}^n }.
\end{array} \right.
\eeq
Note $g_0 = 1$ in the above.
For each given $k$, \reff{min:Rz(f):order=k} is a semidefinite program.
{The length of the moment vector $w$ is $\binom{n+2k}{2k}$.
The vector $\mathscr{V}_{c_i}^{(2k)}[w]$ has length $\binom{n+2k-\deg(c_i)}{2k-\deg(c_i)}$.
The matrix $\mathrm{L}_{c_j}^{(k)}[w]$ has length
$\binom{n+k-\lceil \deg(c_i) \rceil}{k-\lceil \deg(c_i) \rceil}$.
The length of $\mathrm{L}_{\nu, g_j}^{(k,k)}[w]$ is
$\binom{m+k}{k} \cdot \binom{n+k}{k}$.}
The following is the convergence property of
the moment relaxations \reff{min:Rz(f):order=k}.

\begin{theorem}\label{thm:lass}
Let $K$ be  as in \eqref{eq:K}
 {and suppose
the measure $\nu$ satisfies Assumption~\ref{ass:measure}}.
Assume $\qmod{g,\nu} + \ideal{c_{eq}} + \qmod{c_{in}}$ is archimedean.
Then the sequence $(\gamma_k)_k$ of \eqref{min:Rz(f):order=k}
is monotonically increasing and
\[
\gamma_k \, \to \,
f_{\min} \qquad \mbox{as} \qquad k \to \infty.
\]
\end{theorem}
	
\proof{Proof.}
Clearly, the sequence $\gamma_k$ is monotonically increasing
and $\gamma_k \le f_{\min}$ for all $k$.
For all $\eps >0$, the polynomial $f(\x) - f_{\min} + \eps > 0$
on $X \cap K$, so
\[
f(\x) - f_{\min} + \eps \in \QM{g,\nu}_{2k} +
\ideal{c_{eq}}_{2k} + \qmod{c_{in}}_{2k},
\]
for $k$ large enough, by Theorem \ref{thm:univPsatz}.
For each truncated multisequence $w$
that is feasible in \reff{min:Rz(f):order=k}, we have
\[
\mathscr{R}_w( f(\x) - (f_{\min} - \eps ) 1  ) \ge 0.
\]
This implies that
\[
\mathscr{R}_w( f ) \ge  (f_{\min} - \eps )\mathscr{R}_w( 1  )
= f_{\min} - \eps.
\]
So the optimal value $\gamma_k \ge f_{\min} - \eps$.
Since $\eps >0$ can be arbitrarily small,
the limit of $\gamma_k$ must be $f_{\min}$.
\endproof

\subsection{Sampling}
In the expression for the localizing matrix $\mathrm{L}_{\nu,g_j}^{(k,l)}[w]$ in \eqref{min:Rz(f):order=k},
we need the matrix $Y_{ji}^{(l')}$, which then
requires the moments $\int y^\af \,\dd \nu(y)$,
for the chosen  measure $\nu$ with
$\supp{\nu} = Q$. If $Q$ is a
well-known and understood
set (e.g., a box $[-1, 1]^n$, a simplex, a unit ball or a sphere),
the moments can be given by explicit formulas, such as
for the uniformly distributed probability measure.
If $Q$ is not such a convenient set,
the moments $\int_Q y^\af\, \dd \nu(y)$ may not be readily available.
However, this issue can be fixed by sampling.

For a given degree $l$, the moment vector
$\int [y]_{2l}\, \dd \nu(y)$ can always be written as a sample average,
i.e., there exist points $u_1, \ldots, u_N\in Q$ such that
\[
\int [y]_{2l}\, \dd \nu(y) \, = \, \frac{1}{N}
\sum_{i=1}^N [u_i]_{2l} .
\]
This is guaranteed by Caratheodory's theorem \cite[Theorem I.2.3]{Bar02}.
 {To get such sample points $u_i$ can be tricky for some $Q$.
In our computation, we assume they are available from the description of $Q$.
Properties for them to satisfy are discussed in Theorem~\ref{thm:relint}.  }
Interestingly, the above sample average is actually
the moment sequence of a certain measure whose support equals
$\supp\nu=Q$ if the sample size $N$ is large enough.
    {We refer to Remark~\ref{rmk:sample:Q}
   for how large the sample size $N$ should be.}
For a given degree $d$, consider the cone of all possible moment sequences
\beq \label{momcone:Pd}
\mathcal{P}_d  \coloneqq  \Big\{ \int [y]_d\, \dd \mu(y)\, \Big \vert\,
\mu \text{ is a Borel measure on }  \re^m, \;
\supp{\mu} \subseteq Q
\Big\}.
\eeq

Denote the relative interior of $\mathcal{P}_d$
by $\relint ( \mathcal{P}_d )$.

\begin{theorem} \label{thm:relint}
Let $Q$ be a closed set and let $\mathcal{P}_d$ be as above.
For every $\xi \in \relint ( \mathcal{P}_d )$,
there exists a measure $\nu$ on $\re^m$ such that
\beq \label{sample:nu}
\xi \, = \, \int [y]_d \, \dd \nu(y), \quad \supp{\nu} = Q.
\eeq
Moreover, for points $u_1, \ldots, u_D \in Q$,
if $\dim \mathrm{Span}\{ [u_1]_d, \ldots, [u_D]_d \} = \dim \mathcal{P}_d$,
then the sample average
\[
A(u_1, \ldots, u_D) : = \frac{1}{D}( [u_1]_d + \cdots + [u_D]_d )
\]
belongs to the relative interior $\relint  ( \mathcal{P}_d )$.
\end{theorem}

\proof{Proof.}
Consider the subcone
\[
\mathcal{P}'_d  \coloneqq  \Big\{ \int [y]_d\, \dd \mu(y)
\, \Big \vert  \,
\mu \text{ is a Borel measure on } \re^m, \;
\supp{\mu} = Q
\Big\}.
\]
We show that $\mathcal{P}'_d$ is contained in the
relative interior of $\mathcal{P}_d$.
Let $\mathrm{T}$ denote the embedding Euclidean space of
$[y]_d$ over all possible $y \in \re^m$.
Then $\mathcal{P}'_d \subseteq \mathcal{P}_d$ are both convex cones in $\mathrm{T}$.
Let $\ell$ be any linear functional such that
\[
\ell \ge  0 \quad \text{on}\quad  \mathcal{P}_d, \qquad
\ell(\eta) = 0 \quad \text{for some}\quad \eta \in \mathcal{P}'_d .
\]
 {Since $\ell\geq0$ on $\mathcal{P}_d$
and $[y]_d \in \mathcal{P}_d$ for all $y \in Q$, it is evident that
the polynomial $p(y) \coloneqq \ell ([y]_d)$
is nonnegative on $Q$.}
      Let $\mu$ be the Borel measure such that
$\eta = \int [y]_d\, \dd \mu(y)$ and $\supp{\mu} = Q$. Then
\[
0 = \ell(\eta) =  \int \ell( [y]_d )\, \dd \mu(y) = \int p(y)\, \dd \mu(y).
\]
Since $p(y) \ge 0$ on $Q$ and $\supp{\mu} = Q$,
the above implies that $p(y) \equiv 0$ on $Q$,
i.e.,  $\ell \equiv 0$ on $\mathcal{P}_d$.
This shows that every supporting hyperplane of $\mathcal{P}_d$
passing through any point of $\mathcal{P}'_d$
must also contain $\mathcal{P}_d$ entirely.
So $\mathcal{P}'_d$ lies in the relative interior of $\mathcal{P}_d$.
We remark that $\mathcal{P}'_d$ is dense in $\mathcal{P}_d$.
To see this, fix a measure $\nu$ such that $\supp{\nu} = Q$.
Then for every $\xi \in \mathcal{P}_d$ and each integer $k > 0$,
we have $\xi + \frac{1}{k} \int [y]_d\, \dd \nu(y) \in \mathcal{P}'_d$
and it converges to $\xi$ as $k$ goes to infinity.
Since $\mathcal{P}'_d$ is dense in $\mathcal{P}_d$,
they have the same relative interior.

So, every $\xi \in \relint  ( \mathcal{P}_d )$
is the expectation of $[y]_d$
for a certain measure $\nu$ whose support equals $Q$.
Let $\ell$ be a linear functional such that $\ell \ge 0$ on $\mathcal{P}_d$.
If $\ell( A(u_1, \ldots, u_D) ) = 0$, then
\[
\mathrm{Span}\{ [u_1]_d, \ldots, [u_D]_d \}  \subseteq  \ker \ell.
\]
If $\dim \mathrm{Span}\{ [u_1]_d, \ldots, [u_D]_d \} = \dim \mathcal{P}_d$, then
\[
\mathcal{P}_d  \subseteq \ker \ell.
\]
This implies that $\ell \equiv 0$ on $\mathcal{P}_d$.
The above is true for every linear functional $\ell \ge 0$ on $\mathcal{P}_d$.
Therefore, $A(u_1, \ldots, u_D)$ lies in the relative interior of $\mathcal{P}_d$.
\endproof

\begin{remark} \label{rmk:sample:Q}
Note that the measure $\nu$ in \reff{sample:nu}
automatically satisfies the Carleman condition \reff{fix:suppnu=Q} if $Q$ is bounded.
   However, for unbounded $Q$,
we are not sure if \reff{fix:suppnu=Q} still holds.
For the case of unbounded $Q$, we may apply the
homogenization trick to transform to bounded sets.
We refer to the work \cite{HNY01,HNY02} for how to do this.
To summarize, to formulate the localizing matrix
$\mathrm{L}_{\nu, g_j}^{(k,l)}[w]$ in \reff{express:Lhji[z]},
we can select sample points $u_1, \ldots, u_N\in Q$
such that $\mathrm{Span}\{ [u_1]_{2l}, \ldots, [u_N]_{2l} \}$
has maximum dimension, and then let
\beq \label{sample:Yjiell}
Y_{h_{ji}}^{(l')} \, = \, \frac{1}{N}
\sum_{t=1}^N h_{ji}(u_t) [u_t]_{l'} [u_t]_{l'}^\T.
\eeq
\end{remark}

\section{Numerical Experiments}
\label{sc:Num}

This section reports numerical examples to show
the hierarchy of moment relaxations~\reff{min:Rz(f):order=k}
for solving the SIP \reff{primal-SIP}.
The computations are implemented in MATLAB R2023b on a laptop equipped with
a 10th Generation Intel® Core™ i7-10510U processor and 16GB memory.
The moment relaxations are implemented by the software
\texttt{Gloptipoly}~\cite{henrion2009gloptipoly},
which calls the software \texttt{SeDuMi} \cite{strum2001sedumi}
to solve the corresponding semidefinite programs. For the SIP~\reff{primal-SIP},
we use $x^*$ and $f^*$ to denote the global minimizer and the
global minimum value respectively.
The relaxation order is labelled by $k$. For each $k$, we use $w^{(k)}$
to denote the minimizer of \reff{min:Rz(f):order=k}.
The minimum value of \reff{min:Rz(f):order=k} is denoted as $\gamma_k$,
which is a lower bound for the SIP \reff{primal-SIP}.\looseness=-1

The flat extension condition \reff{FEC:w} can also be used to get
minimizers. However, this works only if the moment relaxation
\reff{min:Rz(f):order=k} is tight for solving the SIP.
When \reff{FEC:w} fails, a practical way to get
an approximate minimizer is to let
\[
\hat{x}_k \coloneqq  ( w^{(k)}_{e_1}, \ldots, w^{(k)}_{e_n} ).
\]
The feasibility of the computed point $\hat{x}_k$ is measured as the function value
\[
\dt_{k,j} \, \coloneqq \,
\min\limits_{y\in Q} g_j(\hat{x}_k,y),
\quad j = 1, \ldots, s.
\]
Then $\hat{x}_k$ satisfies the inequality
constraint in \reff{primal-SIP}
if and only if
\[
\dt_k   \coloneqq    \min\limits_{j \in [s] } \dt_{k,j} \, \ge \, 0.
\]
For each example, if the measure $\nu$ is not specified,
we set it to be the normalized Lebesgue measure so that $\nu(Q) = 1$.
The consumed computational time is denoted as {\tt time}.
For neatness of the presentation, all computational results 
are displayed with four decimal digits.

\begin{example}
(i) Consider the following SIP from \cite{coope1985projected,wang2015feasible}:
\beq \label{Moment_1}
\left\{
\begin{array}{cll}
	\min\limits_{x\in \mathbb{R}^2}& \frac{1}{3}x^2_1+x^2_2+\frac{1}{2}x_1\\
	\st  &-\left(1-x_1^2 y^2\right)^2+x_1 y^2+x_2^2-x_2\geq 0\quad\forall y\in Q,
\end{array}
\right.
\eeq
where $Q= [0,1]$. Computational results for Problem~\eqref{Moment_1} are shown in Table~\ref{Tab:Mom1}.
\begin{table}[htb] 
\begin{tabular}{ccccc} \specialrule{.1em}{0em}{0.1em}
$ k $ &  $\hat{x}_k$ &  $\gamma_k$  &  $\dt_k$ & time(s) \\  \hline
$ 3 $ & $(-0.8433, -0.6041)$ & $ 0.1803 $ & $ -0.0310 $ & $ 0.4503 $ \\
$ 4 $ & $(-0.7847, -0.6140)$ & $ 0.1899 $ & $ -0.0090 $ & $ 0.4991 $ \\
$ 5 $ & $(-0.7650, -0.6164)$ & $ 0.1926 $ & $ -0.0036 $ & $ 0.5749 $ \\
$ 6 $ & $(-0.7574, -0.6173)$ & $ 0.1935 $ & $ -0.0017 $ & $ 0.9063 $ \\
$ 7 $ & $(-0.7541, -0.6176)$ & $ 0.1940 $ & $ -9.46\cdot 10^{-4}$ & $ 1.8849 $ \\
\specialrule{.1em}{0em}{0.1em}
\end{tabular}
\caption{Computational results for SIP~\reff{Moment_1}.}
\label{Tab:Mom1}
\end{table}
The true minimizer is $x^* \approx (-0.7500,   -0.6180)$ with
the minimum value $ f^* \approx 0.1945$. %

\smallskip

(ii) Consider the following SIP:
\beq \label{Moment_4}
\left\{
\begin{array}{cll}
	\min\limits_{x\in X}&(x_1-x_2)(x_1-1)+(x_2-x_1)(x_2-1)
	+(x_1-1)(x_2-1)+x^3_1+x^3_2\\
	\st  &x_1x_2y_1y_2-(x_1x_2+x^2_2+0.01)(y_1y_3+y_2+1)-
	x^2_2y_2y_3\geq 0
	\quad \forall y\in Q,
\end{array}
\right.
\eeq
where the sets are
\[
\begin{aligned}
	&X=[-10,10]^2 \,\cap \, \left\{(x_1, x_2) :
	\begin{aligned}
		x_1x_2+x_1+1\geq 0\\
	\end{aligned}
	\right\},\\
	& Q =\left\{ y\in\mathbb{R}^3 :
	y_1\geq 0, \, y_2 \ge 0, \, y_3 \ge 0,  1-y_1-y_2-y_3  \geq 0
	\right\}.
\end{aligned}
\]
As in \cite{lasserre2021simple}, we have the moment formula:
\[
\int_{Q}y^{\alpha_1}_1 y^{\alpha_2}_2 y^{\alpha_3}_3  \, \dd y =
\displaystyle\frac{6\alpha_1!\alpha_2!\alpha_3!}{(|\alpha|+3)!}.
\]
Computational results for Problem~\eqref{Moment_4} are shown in Table~\ref{Tab:Mom4}.
\begin{table}[htb] 
\begin{tabular}{ccccc} \specialrule{.1em}{0em}{0.1em}
$k$ &  $\hat{x}_k$ & $\gamma_k$  &  $\dt_{k}$ & time(s)\\ \hline
$2$ & $( 0.3643,   -0.0327)$ & $ 0.8624$ & $ -0.0017$ & $  0.4783$ \\	
$3$ & $(0.3661,   -0.0340)$ & $ 0.8645$ & $-0.0012$ & $0.6968$  \\
$4$ & $( 0.3671,   -0.0347)$ & $0.8658$ & $-9.28\cdot 10^{-4}$ & $3.1540$ \\
$5$ & $(  0.3677 ,  -0.0351)$ & $ 0.8665$ & $-7.71\cdot 10^{-4}$ & $66.2828$ \\
\specialrule{.1em}{0em}{0.1em}
\end{tabular}
\caption{Computational results for SIP~\reff{Moment_4}.}
\label{Tab:Mom4}
\end{table}
The true minimizer is $x^* \approx (0.3705 ,  -0.0371)$
with the minimum value $f^* \approx 0.8697$.
\end{example}

\begin{example}
(i) Consider the following SIP:
\beq \label{Moment_2}
\left\{
\begin{array}{cll}
	\min\limits_{x\in\mathbb{R}^3}& -x^2_1(100-x_1-x_2)+x^2_2+2x^2_3 \\
	\st  &
	\begin{pmatrix}
		x_1y^2_1-x_1x_2y_1y_2-x_2x_3y^3_2+0.1\\
		x^2_3(y^2_1-y^2_2)+x^2_2y_1y_2+x_1y_2+0.1
	\end{pmatrix}\geq 0\quad\forall y\in Q,\\
\end{array}
\right.
\eeq
where $Q= \{(y_1, y_2) : y^4_1+y^4_2= 1\}$.
We use the sampling as in \reff{sample:Yjiell} to get moments for $\nu$.
Computational results for for Problem~\eqref{Moment_2} are shown in Table~\ref{Tab:Mom2}.
	\begin{table}[htb] 
		\begin{tabular}{ccccr} \specialrule{.1em}{0em}{0.1em}
			$ k $ &  $\hat{x}_k$ & $\gamma_k$ &  $(\dt_{k,1}, \dt_{k,2})$ & time(s) \\	\hline
			$ 2 $ & $(   -0.1253,   -0.0078 ,  -0.0000)$ & $ -1.5728$ & $(-0.0254 ,  -0.0253
			)$ & $ 0.4407$ \\
			$ 3 $ & $(-0.1128 ,  -0.0062,   -0.0000)$ & $  -1.2747
			$ & $( -0.0129 ,  -0.0128)$ &  $0.7351$ \\
			$ 4 $ & $( -0.1016  , -0.0051 ,  -0.0000)$ & $ -1.0340
			$ & $(-0.0016,   -0.0016)
			$ & $4.3025$ \\
			$ 5 $ & $(  -0.1009,   -0.0049,   -0.0000
			)$ & $ -1.0247
			$ & $( -9.52,   -9.03)\cdot 10^{-4} $ & $104.4387$ \\
			\specialrule{.1em}{0em}{0.1em}
		\end{tabular}
		\caption{Computational results for SIP~\reff{Moment_2}.}
		\label{Tab:Mom2}
	\end{table}
The true minimizer is $x^* \approx (-0.1000,  -0.0018, 0.0000)$,
with  minimum value $ f^* \approx -1.0010$. %

\smallskip

(ii) Consider the following SIP:
\beq  \label{Moment_6}
\left\{
\begin{array}{cl}
	\min\limits_{x\in X}& \big( -\sum\limits_{i=1}^4 x^4_i \big) +x^3_1x^2_2+x^2_2x^3_3+
	x_3x^4_4-x^2_1x^2_2+x_1x_2x_3x_4+x_1x_3 \\[4mm]
	\st  &
	y^{\top}\begin{bmatrix}
		x^2_1-x_2x_3 & x_1+x_2x_4 & x^2_3-x_1x_2\\
		x_1+x_2x_4 & x_1-x^2_4 & 1-e^{\top}x\\
		x^2_3-x_1x_2 & 1-e^{\top}x & x_1x_2+x_3x_4
	\end{bmatrix}y\geq 0\quad \forall y \in Q,
\end{array}
\right.
\eeq
where $X=\{x\in\mathbb{R}^4: 4-x^{\top}x\geq 0\}$ and
\[
Q = \left\{ y \in \re^3: \, 1-y^Ty = 0, \, y \ge 0  \right \}.
\]
We apply the sampling as in \reff{sample:Yjiell} to get moments of $\nu$.
Computational results for Problem~\eqref{Moment_6} are shown in Table~\ref{Tab:Mom6}.
\begin{table}[h]
\begin{tabular}{c c c c c c} \specialrule{.1em}{0em}{0.1em}
	$k$ & $\hat{x}_k$ & $\gamma_k$ & $\dt_k$ & time(s)\\	\hline
	$3$ &  $( 1.3903 ,  -0.0000 ,  -0.7310  ,  0.0000
	)$ & $-16.1250$ &$ -3.12 \cdot 10^{-6}$ & $2.0730$ \\
	$4$ &  $(0.1252, 0.0000, -1.9961, 0.0000)$ & $-16.1250$ &
                     $2.62 \cdot 10^{-6}$ &  $192.2331 $\\
	\specialrule{.1em}{0em}{0.1em}
\end{tabular}
\caption{Computational results for SIP~\reff{Moment_6}.}
\label{Tab:Mom6}
\end{table}
The true minimizer is $x^* \approx (0.1252,    0.0000 ,  -1.9961,    0.0000)$,
and the minimum value is $f^* \approx  -16.1250$.
\end{example}

The following are examples where the quantifier set $Q$ is not semialgebraic.
 {For such $Q$, we typically need sampling to get moments of $\nu$.
We refer to Remark~\ref{rmk:sample:Q} for this issue.
Generally, we pick sample points $u_1, \ldots, u_N\in Q$
such that $\mathrm{Span}\{ [u_1]_{2l}, \ldots, [u_N]_{2l} \}$
has maximum dimension.}

\begin{example}
(i) Consider the following SIP
\beq \label{Moment_8}
\left\{
\begin{array}{lll}
	\min\limits_{x\in X}& -x_1x_2x_3+x^3_1+x^2_2+x_3 \\
	\st
	&(x_1x_2+1)y^4_2+(e^{\top}x)y^2_1y_2+(x_1+x_2x_3)y^3_1-0.1\geq 0
	\quad \forall y\in Q,
\end{array}
\right.
\eeq
where the sets
\[
X=\left\{x\in\mathbb{R}^3 \middle\vert\,
\begin{array}{l}
	5-x^{\top}x\geq 0,\\
	x_1x_2-x_3\geq 0
\end{array}
\right \}, \quad
Q = \left\{y\in\mathbb{R}^2\middle\vert\,
\begin{array}{l}
	4-3^{y^2_1}-3^{y^2_2}\geq 0,\\
	3^{y_1}-3^{y_2}-1 \geq 0
\end{array}
\right \}.
\]
We apply the sampling as in \reff{sample:Yjiell} to get moments of $\nu$.
Computational results for Problem~\eqref{Moment_8} are shown in Table~\ref{Tab:Mom8}. 	
	\begin{table}[htb]
		\begin{tabular}{ccccc} \specialrule{.1em}{0em}{0.1em}
			$k$ &  $\hat{x}_k$ & $\gamma_k$ & $\dt_k$ & time(s)\\ \hline
			$2$ & $(  -2.0115,   -0.6861,   -0.6951)$ & $  -7.4037$ & $ -2.7956$ & $0.4233$  \\
			$3$ & $(  0.4350,   -0.3706,   -2.1618)$ & $ -2.2907
			$ & $-7.89\cdot10^{-4}$ & $0.5196$ \\
			$4$ & $(0.4350 ,  -0.3706,   -2.1618
			)$ & $  -2.2907
			$  & $-7.71\cdot 10^{-4}$ & $1.4434$ \\
			$5$ & $( 0.4350,   -0.3707 ,  -2.1618)$& $ -2.2907
			$ &   $ -6.25\cdot 10^{-4}$ & $19.9535$ \\
			\specialrule{.1em}{0em}{0.1em}
		\end{tabular}
		\caption{Computational results for SIP~\reff{Moment_8}.}
		\label{Tab:Mom8}
	\end{table}
The true minimizer is $x^* \approx (0.4353,   -0.3710,   -2.1617)$,
and the minimum value is $f^* \approx  -2.2907$.
They are estimated by applying the $6$th
order Taylor expansion of the exponential function.
			
\smallskip

(ii) Consider the following SIP:
\beq  \label{Moment_9}
\left\{
\begin{array}{lll}
	\min\limits_{x\in X}& x^3_1-x^3_3+x_1x^2_2+(x_2+x^2_3)^2\\
	\st  &
	x_1x_3y_2y_3+x_1x_2y_3+x_2x_3y_1+(x_1+2x_2+x_3)(y_1y_2+2y_3)\geq 0
          \quad
	\forall y\in Q,
\end{array}
\right.
\eeq
where $X=[-1.5, 1.5]^3$ and
\[
Q=\left\{ y\in\mathbb{R}^3 \middle\vert\,
\begin{array}{l}
	2-y^{\top}y\geq 0,\\
	2^{y_3}-2^{y_1}-2^{y_2}\geq 0
\end{array}
\right\}.
\]
We apply the sampling as in \reff{sample:Yjiell} to get moments of $\nu$.
Computational results are shown in Table~\ref{Tab:Mom9}.
\begin{table}[htb]
\begin{tabular}{c c c c c} \specialrule{.1em}{0em}{0.1em}
$k$ & $\hat{x}_k$ & $\gamma_k$ & $\dt_k$ & time(s)\\  \hline
$3$ & $(-1.5000,  1.5000, -0.0000)$ & $-4.5000$ & $-6.76\cdot 10^{-6}$ & $0.7619$ \\
$4$ & $(-1.5000,  1.5000, -0.0000)$ & $-4.5000$ & $-7.99\cdot 10^{-6}$& $ 14.9845$ \\
\specialrule{.1em}{0em}{0.1em}
\end{tabular}
\caption{Computational results for SIP~\reff{Moment_9}.}
\label{Tab:Mom9}
\end{table}
The true minimizer is $x^* \approx (-1.5000,    1.5000 ,   0.0000)$
and $f^* \approx   -4.5000$.	
They are estimated by applying the $6$th order Taylor expansion
of exponential functions.

\smallskip

{(iii) Consider the following SIP:
\beq  \label{Moment_Zm=1}
\left\{
\begin{array}{rl}
	\min  & x_1^3+x_2^3 \\
	\st  &   4y^4 -1 - y^2(2y^2-1) x_1^2 - y^2(2y^2+1) x_2^2   \geq 0
   \qquad   \forall y\in Q,
\end{array}
\right.
\eeq
where $Q = {\Z_+}$.
We select the same measure $\nu$ as in Example~\ref{exm3.10:Q=Z}.
The computational results are shown in Table~\ref{Tab:Ex_6_3_iii}. 	
For this SIP, the true minimum value $f^* \approx -2.8284$
and the true minimizer $x^* \approx (-1.4142, 0.0000)$.
We remark that there are numerical issues for solving
the moment relaxation~\reff{min:Rz(f):order=k}
when the relaxation order $k \ge 7$.
\begin{table}[htb]
\begin{tabular}{cccc} \specialrule{.1em}{0em}{0.1em}
$k$ &  $\hat{x}_k$ & $\gamma_k$   & time(s)\\ \hline
$2$ & $(-1.4561, -0.0000)$ & $ -3.0874 $  &  $0.7859$  \\
$3$ & $(-1.4322, -0.0000)$ & $ -2.9375 $  &  $0.9379$  \\
$4$ & $(-1.4246, -0.0000)$ & $  -2.8915 $  &  $0.8905$  \\
$5$ & $(-1.4211, -0.0000)$ & $ -2.8701 $  &  $1.0247$  \\
$6$ & $(-1.4191, -0.0000)$ & $ -2.8581 $  &  $1.4835$  \\
$7$ & $(-1.4179, -0.0002)$ & $ -2.8506 $  &  $2.4642$  \\
$8$ & $(-1.4156, -0.0002)$ & $  -2.8424 $  &  $5.8728$  \\
	\specialrule{.1em}{0em}{0.1em}
\end{tabular}
\caption{Computational results for the SIP~\reff{Moment_Zm=1}.}
\label{Tab:Ex_6_3_iii}
\end{table}
      }
\end{example}

The following are examples where
the quantifier set $Q$ is a union of several closed sets, {say,
\[ Q \, = \, Q_1 \cup \cdots \cup Q_l,  \quad \text{for each} \quad Q_i \subseteq \re^m.   \]
For each $i$, let $\nu$ be a measure on $\re^m$ such that $\supp{\nu_i} = Q_i$.
Then $\nu \coloneqq \nu_1 + \cdots + \nu_l$ is a measure such that $\supp{\nu} = Q$.
By the definition, one can see that
\[ \qmod{g,\nu} \, = \,\qmod{g,\nu_1} + \cdots + \qmod{g,\nu_l} .    \] }

\begin{example}
(i) Consider the following SIP:
\beq  \label{Moment_11}
\left\{
\begin{array}{cll}
	\min\limits_{x\in\mathbb{R}^3}& (x^2_1+1.8x^2_3)^2+x_1x_2x_3+x^3_1-2x^3_2-4x_3\\
	\st  & \begin{pmatrix}
		x_1x_2y_1y_2-x_2x_3(y_1+y^2_3)-0.01\\
		x^2_3y^2_1-x^2_2y_2y_3+x^2_1(y_1+y_3-0.1)
	\end{pmatrix}\geq 0\quad\forall y\in Q,\\
\end{array}
\right.
\eeq
where $Q = Q_1 \cup Q_2$ is the union of the following two sets:
\[
\begin{array}{rcl}
	Q_1  &=&  \left\{y\in\mathbb{R}^3 :
	\begin{aligned}
		(y_1- 1)^2+(y_2-1)^2+(y_3-1)^2\leq 1
	\end{aligned}\right\},  \\
	Q_2 &=& \{ y\in\mathbb{R}^3 :
	\begin{aligned}
		(y_1- 1)^2+y^2_2+(y_3-1)^2\leq 1
	\end{aligned} \}.
\end{array}
\]
We apply the sampling as in \reff{sample:Yjiell} to get moments of $\nu$.
Computational results are shown in Table~\ref{Tab:Mom11}.
\begin{table}[htb]
\begin{tabular}{ccccr} \specialrule{.1em}{0em}{0.1em}
	$k$ & $\hat{x}_k$ & $\gamma_k$ & $(\dt_{k,1}, \dt_{k,2} )$ &  time(s) \\ \hline
	$2$ & $  (  1.8793  ,  2.2691,   -0.2007
	)$ & $ -3.7910 $ & $(  -4.4703  , -6.6000
	)$ & $0.5968$\\
	$3$ & $(  0.0047,   -0.0231,    0.6758
	)$ & $ -2.0274$ & $(-4.44,-0.58)\cdot 10^{-3}$ & $ 1.1363$\\
	$4$ & $( 0.0047 ,  -0.0230 ,   0.6758
	)$ & $ -2.0274
	$ & $(-4.42,-0.58)\cdot 10^{-3}$ & $   48.7642
	$\\
	\specialrule{.1em}{0em}{0.1em}
\end{tabular}
\caption{Computational results for SIP~\reff{Moment_11}.}
\label{Tab:Mom11}
\end{table}

(ii) Consider the following SIP:
\beq  \label{Moment_12}  
\left\{
\begin{array}{cll}
	\min\limits_{x \in X} & (x^2_1-x_2)^2-3x_1x^2_2+3x^3_1\\
	\st  & 	-x_1x_2(y^2_1+2y^2_3)+x^2_2(y_1-y_2y_3)+ \\
         &  \qquad \qquad  2y_1y_3-e^T x-1.4\geq 0  \quad\forall y\in Q,
\end{array}
\right.
\eeq
where $X=\{x\in\mathbb{R}^2 : 8- x^T x \geq 0\}$ and
\[
Q= \left\{y\in\mathbb{R}^3 \middle\vert\,
\begin{array}{l}
	10-y^T y\geq 0,\\
	|y_1|+|y_2|+|y_3|-1\geq 0
\end{array}
\right\} .
\]
Note that $Q$ is a union of $8$ basic closed semialgebraic sets, that is,
\[
Q \quad  = \bigcup_{ \substack{ s_1, s_2, s_3 \in \{ -1, 1\}  }  }
Q_{s_1, s_2, s_3 } \coloneqq \left\{y\in\mathbb{R}^3 \middle\vert\,
\begin{array}{l}
	10-y^T y\geq 0, \\
	s_1 y_1 \ge 0, s_2 y_2 \ge 0, s_3 y_3 \ge 0, \\
	s_1 y_1 + s_2 y_2 + s_3 y_3 -1  \ge  0
\end{array}
\right\} .
\]
We apply the sampling as in \reff{sample:Yjiell} to get moments of $\nu$.
Computational results are shown in Table~\ref{Tab:Mom12}.
\begin{table}[htb]
\begin{tabular}{c c c c c}  \specialrule{.1em}{0em}{0.1em}
$k$ &  $\hat{x}_k$ & $\gamma_k$ & $\dt_k$ &  time(s)\\  \hline
$2$ & $(-2.4791 ,   0.7284)$ & $-12.4135$ & $-6.46 \cdot 10^{-4}$ & $0.4755$ \\
$3$& $(-2.4791 ,   0.7284)$ & $-12.4135$ & $-6.33 \cdot 10^{-4}$ & $ 0.7434$ \\
$4$& $(-2.4791,    0.7284)$ & $ -12.4135$ &  $-6.36 \cdot 10^{-4}$ & $ 2.8389$ \\
$5$& $(   -2.4779,    0.7272 )$ & $  -12.4099$ &  $ 4.01 \cdot 10^{-4}$ & $53.6131$ \\
\specialrule{.1em}{0em}{0.1em}
\end{tabular}
\caption{Computational results for SIP~\reff{Moment_12}.}
\label{Tab:Mom12}
\end{table}
\end{example}

We remark that the flat extension condition \reff{FEC:w}
can also be used to get minimizers when it holds.
This happens only if the moment relaxation is tight
for solving the SIP. See the following example.

\begin{example}\label{tms_2}
Consider the following SIP:
\beq  \label{multi-minimizer}
\left\{
\begin{array}{cll}
	\min\limits_{x\in X}& -x^2_1x^2_2\\
	\st  & (x_1+x_2)y^2_2-x_1x_2(y_1y_2+1)\geq 0\quad\forall y\in Q,\\
\end{array}
\right.
\eeq
where $X=\{x\in\mathbb{R}^2: 1-x^Tx\geq 0\}$ and \[Q=\{y\in\mathbb{R}^2: |y_1|+|y_2|\leq 1\}.\]
For the relaxation order $k=2$, we get the optimal $w^*$ such that
\[
\mathrm{H}^{(2)}[w^*] \, = \,
\begin{bmatrix}
1 & 0 & 0 & \frac{1}{2} & -\frac{1}{2} & \frac{1}{2} \\[3pt]
0 & \frac{1}{2} & -\frac{1}{2} & 0 & 0 & 0 \\[3pt]
0 & -\frac{1}{2} & \frac{1}{2} & 0 & 0 & 0 \\[3pt]
\frac{1}{2} & 0 & 0 & \frac{1}{4} & -\frac{1}{4} & \frac{1}{4} \\[3pt]
-\frac{1}{2} & 0 & 0 & -\frac{1}{4} & \frac{1}{4} & -\frac{1}{4} \\[3pt]
\frac{1}{2} & 0 & 0 & \frac{1}{4} & -\frac{1}{4} & \frac{1}{4}
\end{bmatrix} .
\]
The flat extension \reff{FEC:w} holds. Indeed,
we can get $w=\frac{1}{2}([u^*_1]_4+[u^*_2]_4)$ for points
\[
u^*_1=\Big(-\frac{1}{\sqrt{2}},\frac{1}{\sqrt{2}}\Big), \quad 
u^*_2=\Big(\frac{1}{\sqrt{2}},-\frac{1}{\sqrt{2}}\Big).
\]
They are both minimizers for this SIP.
\end{example}

\section{Conclusions and Discussions}
\label{sec:con}

We study Positivstellens\"atze and moment problems for sets
that are given by universal quantifiers.
For the set $K$ as in \reff{eq:K} given by a universal quantifier
$y \in Q$, we discuss representation of polynomials
that are positive on $K$. Let $\nu$ be a measure satisfying the Carleman
condition~\reff{ass:measure}. When the quadratic module $\QM{g,\nu}$ is archimedean,
we show in Theorem~\ref{thm:univPsatz} that a polynomial $f(\x)$  positive on $K$
 {as in \reff{eq:K}} must be  in $\QM{g,\nu}$.
For the non-archimedean case, we give a similar result in Corollary~\ref{cor:perturb}.
We also study $K$-moment problems for the set $K$.
Necessary and sufficient conditions
for a full (or truncated) multisequence to admit a representing measure
supported in $K$ are given. In particular, the classical flat extension
theorem is generalized for truncated moment problems with such a set $K$.
These results are presented in Theorems~\ref{theorem:mp}, \ref{thm:shift}
and \ref{thm:truncMP}, respectively.
These new Positivstellens\"atze and moment problems
can be applied to solve semi-infinite optimization (SIP).
For the SIP~\reff{primal-SIP}, a hierarchy of moment relaxations
\reff{min:Rz(f):order=k} is proposed to solve it.
Its convergence is shown in Theorem~\ref{thm:lass}.
Various examples for semi-infinite optimization
are demonstrated in Section~\ref{sc:Num}.

Our work leads to many intriguing questions to explore in the future.
For instance, without assuming archimedeanness,
is there a preordering version of Theorem \ref{thm:univPsatz}? Equivalently,
does there exist a clean algebraic reformulation of
the compactness (or emptiness) of the set $K$ given with a universal quantifier as in \eqref{eq:K}?
Is  there an analog of the Krivine-Stengle Positivstellensatz for such sets $K$?
It would also be interesting to establish the universal Positivstellens\"atze
for matrix-valued polynomials and matrix-valued constraints.
In Theorem~\ref{thm:truncMP}, the condition~\reff{Lgj[z]>=0:k=1:d}
is assumed to hold for all $l=1, 2, \ldots$. If it holds
for only finitely many $l$, the conclusion
of Theorem \ref{thm:truncMP} may not hold.
It would be interesting to find a \emph{finite} set of conditions
for a truncated multisequence to
admit a representing measure supported in $K$.
Finally, in their previous joint work,
the second and third author \cite{KN20} gave Positivstellens\"atze
and solvability criteria for moment problems for sets given with \emph{existential} quantifiers.
A major future task will be to give a common extension of the results from \cite{KN20}
and the present paper, that is, Positivstellens\"atze and  moment problems
for sets given with a combination of universal and existential quantifiers.

\bigskip \noindent
{\bf Acknowledgment.}
 {We thank the anonymous referees for their valuable comments and suggestions.}
Xiaomeng Hu and Jiawang Nie are partially supported by the NSF grant DMS-2110780.
Igor Klep is supported by
the Slovenian Research Agency
program P1-0222 and grants J1-50002, J1-2453, N1-0217,
J1-3004, and was
partially supported by the Marsden Fund Council
of the Royal Society of New Zealand.
\CRed{Igor's work was partly performed within the project COMPUTE,
funded within the QuantERA II Programme that has received funding from the EU’s H2020 research
and innovation programme under the GA No 101017733
\raisebox{-.5mm}{\worldflag[width=4mm]{EU}}.}
\looseness=-1

\ifx\foo\undefined

\end{document}